\documentclass[11pt]{article}

\usepackage{amssymb,amsmath}
\usepackage{graphicx}
\usepackage{slashbox}
\usepackage{arydshln}
\usepackage{rotating}
\usepackage{natbib}

{\catcode `\@=11 \global\let\AddToReset=\@addtoreset}
\AddToReset{equation}{section}

\newtheorem{cor}{\bf Corollary}[section]
\newtheorem{prop}{\bf Proposition}[section]
\newtheorem{defn}{\bf Definition}[section]

\newtheorem{@rem}{\bf Remark}[section]
\newenvironment{rem}{\begin{@rem}\rm}{\end{@rem}}
\newtheorem{@ex}{\bf Example}[section]
\newenvironment{ex}{\begin{@ex}\rm}{\end{@ex}}

\def\E{\mathbb{E}}

 \def\d{\ \mathrm{d}}

\def\phi{\varphi}

\newcommand\noi{\noindent}

\def\mathsf{\bf}

\def\R{\mathbb{R}}
\def\Z{\mathbb{Z}}

\def\i{\mathrm i}
\def\d{\mathrm d}
\def\e{\mathrm e}

\def\E{\mathrm E}
\def\P{\mathrm P}

\def\text{\mbox}
\def\veps{\varepsilon}
\def\1{{\bf 1}}

\newcommand\beqn{\begin{displaymath}}  
\newcommand\eeqn{\end{displaymath}}

\begin{document}

\title{Detection of non-constant long memory parameter\thanks{The second and fourth authors
are supported by a grant (No.\ MIP-11155) from the Research Council of Lithuania.}}

\author{Fr\'{e}d\'{e}ric Lavancier$^{1}$ \ (frederic.lavancier@univ-nantes.fr,\smallskip\\
Remigijus Leipus$^{2,3}$ \ (remigijus.leipus@mif.vu.lt),\smallskip\\
  Anne Philippe$^{1}$ \ (anne.philippe@univ-nantes.fr),\smallskip\\ 
  Donatas Surgailis$^{2,3}$  \ (donatas.surgailis@mii.vu.lt)
\medskip\\
\\{\small\it $^{1}$Laboratoire de Math\'{e}matiques Jean Leray, Universit\'{e} de Nantes, France,}\\
 {\small\it $^{2}$Faculty of Mathematics and
Informatics, Vilnius University, Lithuania,}\\
{\small\it $^{3}$Institute of Mathematics and Informatics, Vilnius
University, Lithuania} }

\date{\empty}

\maketitle


 {\noindent {\bf Abstract}\ \
This article deals with detection of nonconstant long memory parameter in time series. The null
hypothesis presumes stationary or nonstationary time series with constant long memory parameter,
typically an $I(d)$ series with $d>-.5$. The alternative corresponds to an increase in persistence
and includes in particular an abrupt or gradual change from $I(d_1)$ to $I(d_2)$, $-.5<d_1<d_2$.
We discuss several test statistics based on the ratio of forward and backward sample
variances of the partial sums. The consistency of the tests is proved under a very general setting.
We also study the behavior of these test statistics for some models with changing memory parameter.
A simulation study shows that our testing procedures have good finite sample properties and turn out to be
more powerful than the KPSS-based tests considered in some previous works.
\bigskip


\bigskip

\section{Introduction}


The present paper discusses statistical tests for detection of non-constant memory parameter of
time series versus the null hypothesis that this parameter has not changed over time.
As a particular case, our framework includes testing the null hypothesis that the observed series is
$I(d)$ with constant $d>-.5$, against the alternative hypothesis that $d$ has changed, together with a
rigorous formulation of the last change. This kind of testing procedure is the basis to study the dynamics
of persistence, which is a major question in economy (see \cite{kumar}, \cite{hassler3}, \cite{Kruse:2008}).

In a parametric setting and for stationary series ($|d|<.5$), the problem of testing for a single change
of $d$ was first investigated by \cite{beran}, \cite{horvath:shao:1999}, \cite{horvath}, \cite{yamaguchi}
(see also \cite{lavielle}, \cite{Kokoszka:Leipus-01}).
Typically, the sample is partitioned into two parts and $d$ is estimated on each part.
The test statistic is obtained by maximizing the difference of these estimates over all such partitions.
A similar approach for detecting multiple changes of $d$ was used in \cite{shimotsu}
and \cite{bardet} in a more general semiparametric context.

The above approach for testing against changes of $d$ appears rather natural although
applies to abrupt changes only and involves (multiple) estimation of $d$ which is not
very accurate if the number of observations between two change-points is not large enough; moreover, estimates of $d$ involve
bandwidth or some other tuning parameters and are rather sensitive to the short memory spectrum of the process.


On the other hand, some regression-based Lagrange Multiplier procedures have been recently discussed in
\cite{hassler2} and \cite{martins}. The series is first filtered by $(1-L)^d$, where $L$ is the
lag operator and $d$ is the long memory parameter under the null hypothesis, then the resulting series is subjected to a (augmented) Lagrange Multiplier test for fractional integration, following the pioneer works by Robinson (\citeyear{robinson91}, \citeyear{robinson94}).
The filtering step can be done only approximatively and involves in practice an estimation of $d$. This is certainly the main reason for the size distortion that can be noticed in the simulation study displayed in \cite{martins}.

In a nonparametric set up, following \cite{Kim:2000}, \cite{Kim_B-F_A:2002} proposed several tests (hereafter referred to as
Kim's tests), based on the ratio
\begin{equation} \label{kim0}
 {\cal K}_n(\tau) :=  \frac{U^*_{n-\lfloor n\tau\rfloor}(X)}{U_{\lfloor n\tau\rfloor}(X)}, \qquad \tau \in [0,1],
\end{equation}
where
\begin{eqnarray}\label{kim00}
&U_k(X) := \frac{1}{k^2} \sum_{j=1}^k \big(S_j - \frac{j}{k}S_k\big)^2, \qquad U^*_{n-k}(X) :=  \frac{1}{(n-k)^2} \sum_{j=k+1}^n \big(S^*_{n-j+1} - \frac{n-j+1}{n-k}S^*_{n-k}\big)^2
\end{eqnarray}
are estimates of the second moment of forward and backward de-meaned partial sums
$$
\frac{1}{k^{1/2}}\Big(S_j - \frac{j}{k}S_k\Big), \, j=1,\dots, k \qquad  \text{and}  \qquad
\frac{1}{(n-k)^{1/2}}\Big(S^*_{n-j +1} - \frac{n-j+1}{n-k}S^*_{n-k}\Big), \, j=k+1, \dots, n,
$$
on intervals $[1, 2, \dots, k] $ and $[k+1, \dots, n]$, respectively. Here and below, given a sample $X=(X_1, \dots, X_n)$,
$$
S_{k} := \sum_{j=1}^{k} X_j, \qquad  S^*_{n-k} := \sum_{j=k+1}^{n} X_j
$$
denote the forward and backward partial sums processes.
Originally developed to test for a change from $I(0)$ to $I(1)$
(see also \cite{Busetti_Taylor:2004}, \cite{Kim_B-F_A:2002}), Kim's statistics were extended
in \cite{Hassler_Scheithauer:2011}
to detect a change from $I(0)$ to $I(d)$, $d>0$.
A related, though different approach based on the so-called CUSUM statistics, was used in \cite{Leybourne_Taylor_Kim:2007}
and \cite{Sibbertsen_Kruse:2009} to test for a change from stationarity ($d_1<.5$) to
nonstationarity ($d_2 >.5$), or vice versa.

The present work extends Kim's approach to detect an abrupt or gradual change
from $I(d_1)$ to $I(d_2)$, for any
$-.5 < d_1 < d_2$ with exception of values $d_1, d_2 \in \{.5, 1.5,\dots\}$ (see Remark \ref{rem3.2} for an explanation
of the last restriction).
This includes both stationary and nonstationary null (no-change) hypothesis which is important for
applications since nonstationary time series with $d >.5$ are common in economics.
 Although our asymptotic results (Propositions~\ref{consist}, \ref{general} and Corollary~\ref{consist0})
 are valid for the original Kim's statistics, see Remark~\ref{kpss},
 we modify Kim's ratio (\ref{kim0}), by replacing the second sample moments
 $U_k(X), \, U^*_{n-k}(X)$ in (\ref{kim00}) of backward and forward partial sums
by the corresponding empirical  variances $V_k(X), \, V^*_{n-k}(X)$ defined at
(\ref{varS}) below.
This modification is similar to the difference between the KPSS and the V/S tests, see \cite{GKLT-JE00}, and leads
to a more powerful testing procedure (see Tables~\ref{compa500}--\ref{compa5000}). It is important to note that the ratio-based statistics
 discussed in our paper, as well as the original Kim's statistics, do not require an estimate of $d$ and do not depend on any tuning parameter apart from the choice of the testing interval ${\cal T} \subset (0,1)$. However, the limiting law under the null hypothesis depends on $d$, hence the computation of the quantile defining the critical region requires a weakly consistent estimate of the
memory parameter $d$.

The paper is organized as follows. Section~\ref{sec:hyp} contains formulations
of the null and alternative hypotheses, in terms of joint convergence of forward
and backward partial sums processes, and describes a class of $I(d)$ processes
which satisfy the null hypothesis. Section~\ref{sec:testing-procedure} introduces the ratio
statistics $W_n, I_n$ and $R_n$ and derives their limit distribution under the
null hypothesis. Section~\ref{sec:consistency} displays  theoretical results, from
which the consistency of our testing procedures is derived.  Section~\ref{sec:appl-fract-integr} discusses the behavior of our
statistics under alternative hypothesis.
Some fractionally integrated models with constant or changing memory parameter are considered and
the behavior of the above statistics for such models is studied. Section~\ref{sec:lin-trend}
extends the tests of Section~\ref{sec:testing-procedure} to
the case when observations contain a linear trend. Section~\ref{sec:simulation}
contains simulations of empirical size and power of our testing procedures.
All proofs are collected in Section~\ref{sec:app}.

\section{The null and alternative hypotheses}\label{sec:hyp}

Let $X=(X_1,\dots,X_n)$ be a sample from a time series $\{X_j\} =
\{X_j, j = 1,2,\dots\}$.
Additional assumptions about $\{X_j\}$ will be specified later.
Recall the definition of forward and backward partial sums
processes of $X$:
$$
S_k=S_{k}(X)= \sum_{j=1}^{k}X_j, \qquad S^*_{n-k}=S^*_{n-k}(X)=
\sum_{j=k+1}^{n}X_j.
$$
Note that backward sums can be expressed via forward sums, and
vice versa: $S^*_{n-k} = S_n - S_k$, $S_k = S^*_n - S^*_{n-k}$.

For $0\le a<b\le 1$, let us denote by  $D[a,b]$ the Skorokhod
space of all cadlag (i.e.\ right-continuous with left limits) real-valued functions
defined on interval $[a,b]$.
 In this article, the space  $D[a,b]$  and the product space
$D[a_1,b_1]\times D[a_2,b_2]$, for any $0\leq a_i<b_i\leq 1$,
$i=1,2$, are all endowed with the uniform topology and the
$\sigma$-field generated by the open balls (see \cite{pollard}).
The weak convergence of random elements in such spaces is denoted
'$\longrightarrow_{D[a,b]}$' and '$\longrightarrow_{D[a_1,b_1]\times
D[a_2, b_2]}$', respectively; the weak convergence of finite-dimensional distributions
is denoted '$\longrightarrow_{fdd}$';
the convergence in law and in probability of
random variables are denoted
'$\longrightarrow_{\rm law}$' and '$\longrightarrow_p$', respectively.

The following hypotheses are clear particular cases of our more general
hypotheses ${\bf H_0}$, ${\bf H_1}$ specified later. The null hypothesis
below involves the classical type I fractional Brownian motion in the limit
behavior of the partial sums, which is typical for linear models with long memory.
Recall that a type I fractional Brownian motion $B^{\rm I}_{d+.5}=\{B^{\rm I}_{d+.5}(\tau), \tau\ge 0\}$ 
with Hurst parameter $H=d +.5 \in (0,2)$, $H\ne 1$ is defined by
\begin{eqnarray}\label{fbmI}
 B^{\rm I}_{d+.5}(\tau)&:=&\begin{cases}
\frac{1}{\Gamma(d+1)}\int_{-\infty}^\tau \big((\tau-u)^d - (-u)^d_+\big)\d B(u),  &-.5 < d < .5, \\
\int_0^\tau B^{\rm I}_{d-.5}(u) {\d}u, & .5<d<1 .5, 
\end{cases}
\end{eqnarray}
where $(-u)_+ := (-u)\vee 0$ and $\{B(u), u \in \R\}$ is a standard Brownian motion with
zero mean and variance $\E B^2(u) = |u|$. Let $\lfloor x\rfloor $ denote the integer part of the real number
$x \in \R$.

\bigskip

\noindent ${\bf H_0[I]}$: There exist $d \in (-.5,1.5)$, $d\ne .5$,
$\kappa>0$
and a normalization $A_{n}$ such that
\begin{eqnarray}\label{Dconv00}
 n^{-d -.5} \big(S_{\lfloor n\tau\rfloor} - \lfloor n\tau\rfloor A_{n} \big)
&\longrightarrow_{D[0,1]}& \kappa B^{\rm I}_{d+.5} (\tau),\quad n\to\infty.
\end{eqnarray}

 \bigskip

\noindent ${\bf H_1[I]}$: There exist $0 \le \upsilon_0 <  \upsilon_1 \le
1$, $d > -.5$, and a normalization $A_n$ such that
\begin{eqnarray}\label{Dconv11}
 \Big(n^{-d-.5} \big(S_{\lfloor n\tau_1\rfloor}  - \lfloor n\tau_1\rfloor A_n\big), \, n^{-d-.5} \big(S^*_{\lfloor n\tau_2\rfloor} -
\lfloor n\tau_2\rfloor A_n\big)
\Big)
&\longrightarrow_{D[0, \upsilon_1] \times D[0,
1-\upsilon_0]}&\big(0, Z_2(\tau_2)\big),
\end{eqnarray}
as $n\to\infty$, where
$\{Z_2(\tau), \tau\in [1-\upsilon_1,1-\upsilon_0]\}$ is a
nondegenerate a.s.\ continuous Gaussian process.
\bigskip

Here and hereafter, a random element $Z$ of $D[a,b]$ is called {\it nondegenerate}
if it is not identically zero on the interval
$[a,b]$ with positive probability, in other words, if
$\P(Z(u) = 0, \forall u \in [a,b]) = 0$.

Typically, the null hypothesis ${\bf H_0[I]}$ is satisfied by  $I(d)$
series (see Definition~\ref{defId}).
In Section~\ref{sec:type-i-fractional}
we give a general family of linear processes satisfying ${\bf H_0[I]}$
including stationary and nonstationary processes. See also \cite{taqqu}, \cite{giraitisAAP} and the review paper
\cite{giraitisHandbook} for some classes of non-linear stationary processes (subordinated Gaussian processes and stochastic volatility
models) which satisfy ${\bf H_0[I]}$ for $ 0< d < .5$.
The alternative hypothesis corresponds to the processes changing from $I(d_1)$
to $I(d_2)$ processes (see Section~\ref{sec:fract-integr-models} for examples).
\medskip

Let us  give a first example based on the well-known FARIMA model.

\begin{ex}\label{ex1}
{\rm A FARIMA($0,d,0$) process $\varepsilon_t(d) = \sum_{s=0}^\infty \pi_s(d) \zeta_{t-s}$ with $-.5<d<.5$ satisfies
assumption ${\bf H_0[I]}$ with $\kappa =1, \,A_n =0$.
 Here, $\pi_s(d), s=0,1,\dots$ are the moving-average coefficients (see \eqref{farima})
and $\{\zeta_t\}$ is a Gaussian white noise with zero mean
and unit variance.
Moreover, for two different memory parameters $-.5< d_1 < d_2 < .5$, we can construct a process satisfying ${\bf H_1[I]}$ by
\begin{equation}\label{changed}
 X_t  := \begin{cases} \varepsilon_t (d_1),
   & \text{$ t\leq \lfloor n\theta^*\rfloor, $}\\
\varepsilon_t(d_2),
&  \text{$ t > \lfloor n\theta^*\rfloor,$}
\end{cases}
\end{equation}
where $\theta^* \in (0,1)$. The process in (\ref{changed}) satisfies (\ref{Dconv11}) with
$d = d_2, \, A_n = 0,  \, \upsilon_0 = 0, \, \upsilon_1 = \theta^*,$ and
$Z_2(\tau) =  B^{\rm I}_{d+.5} (1) -   B^{\rm I}_{d+.5} (\theta^* \vee (1-\tau)), \, \tau \in [0,1]$.
In the case of $.5<d<1.5$, FARIMA(0,$d$,0) process in \eqref{changed} is defined by $\varepsilon_t(d)=\sum_{i=1}^t Y_i$,
where $\{Y_t\}$ is a stationary FARIMA(0,$d-1$,0).}
\end{ex}

The testing procedures of Section~\ref{sec:testing-procedure} for testing
the hypotheses  ${\bf H_0[I]}$ and ${\bf H_1[I]}$ can be extended to
more general context. We formulate these `extended' hypotheses as follows.

\medskip

\noindent ${\bf H_0}$: There exist normalizations $\gamma_{n} \to \infty $ and $A_{n}$ such that
\begin{eqnarray}\label{Dconv1}
 \gamma^{-1}_{n} \big(S_{\lfloor n\tau\rfloor } -  \lfloor n\tau\rfloor A_{n}\big)
 &\longrightarrow_{D[0,1]}& Z(\tau),
\end{eqnarray}
where $\{Z(\tau), \tau \in [0,1]\}$ is a nondegenerate a.s.\
continuous random process.

\bigskip

\noindent ${\bf H_1}$: There exist $0 \le \upsilon_0 <  \upsilon_1 \le
1$ and normalizations $\gamma_{n} \to \infty$ and $A_n$  such that
\begin{eqnarray}\label{Dconv}
\Big(\gamma^{-1}_{n} \big(S_{\lfloor n\tau_1\rfloor}  -  \lfloor n\tau_1\rfloor A_n\big), \, \gamma^{-1}_{n} \big(S^*_{\lfloor n\tau_2\rfloor} -
 \lfloor n\tau_2\rfloor A_n\big)
\Big)
&\longrightarrow_{D[0, \upsilon_1] \times D[0,
1-\upsilon_0]}&\big(0, Z_2(\tau_2)\big),
\end{eqnarray}
where $\{Z_2(\tau), \tau \in [1-\upsilon_1,1-\upsilon_0]\}$ is a
nondegenerate  a.s.\ continuous random process.
\bigskip

Typically, normalization $A_{n} = \E X_0$ accounts for centering of observations
and does not depend on $n$. Assumptions ${\bf H_0}$ and ${\bf H_1}$ represent
very general forms of the null (`no change in persistence of $X$') and
the alternative (`an increase in persistence of $X$') hypotheses. Indeed, an
increase in persistence of $X$ at time $k_* = \lfloor n \upsilon_1\rfloor$
typically means that   
forward partial sums $S_j, j \le k_*$ grow at a slower rate $\gamma_{n1}$ compared with
the rate of growth $\gamma_{n2}$ of backward sums $S^*_j, j \le n-k_* $. Therefore, the former sums
tend to a degenerated process $Z_1(\tau) \equiv 0$, $\tau \in [0,\upsilon_1]$ under the
normalization $\gamma_n = \gamma_{n2}$. Clearly, ${\bf H_0}$ and ${\bf H_1}$ are not limited to stationary
processes and allow infinite variance processes as well. While these assumptions are sufficient for
derivation of the asymptotic distribution and consistency of our tests, they need to be specified
in order to be practically implemented. The hypothesis ${\bf H_0[I]}$ presented before is one
example of such specification and involves the type I fBm. Another example involving the type II
fBm is presented in Section~\ref{sec:typeII}.

\section{The testing procedure}\label{sec:testing-procedure}

\subsection{The test statistics}\label{sec:test-statistics}

Analogously to (\ref{kim0})--(\ref{kim00}), introduce the corresponding partial sums' variance estimates
\begin{eqnarray}\label{varS}
V_k(X)& :=&  \frac{1}{k^2} \sum_{j=1}^k \Big(S_j -
 \frac{j}{k}S_k\Big)^2 - \bigg(\frac{1}{k^{3/2}} \sum_{j=1}^k
 \Big(S_j - \frac{j}{k}S_k\Big)\bigg)^2,\nonumber\\
V^*_{n-k}(X) & :=& \frac{1}{(n-k)^2} \sum_{j=k+1}^n \Big(S^*_{n-j+1} - \frac{n-j+1}{n-k}\,S^*_{n-k}\Big)^2  \\
 &&\hspace {.2cm} - \
\bigg(\frac{1}{(n-k)^{3/2}} \sum_{j=k+1}^n \Big(S^*_{n-j+1} - \frac{n-j+1}{n-k}\,S^*_{n-k}\Big)\bigg)^2 \nonumber
\end{eqnarray}
and the corresponding `backward/forward variance ratio':
\begin{equation} \label{calL}
{\cal L}_n(\tau)\ :=\  \frac{V^*_{n-\lfloor n\tau \rfloor}(X)}{V_{\lfloor n\tau \rfloor}(X)}, \qquad \tau \in [0,1].
\end{equation}
For a given {\it testing  interval} ${\cal T} = [\underline\tau,\overline\tau] \subset (0,1)$, define
the analogs of the `supremum' and `integral' statistics of \cite{Kim:2000}:
\begin{eqnarray}\label{memstat}
W_n(X)\ :=\  \sup_{\tau\in {\cal T}} {\cal L}_n(\tau), \qquad  \  I_n(X) \ := \ \int_{\tau\in {\cal T}} {\cal L}_n(\tau) {\d} \tau.
\end{eqnarray}
We also define the analog of the ratio statistic introduced in \cite{Sibbertsen_Kruse:2009}:
\begin{eqnarray}\label{RX}
R_n(X)&:=&\frac{\inf_{\tau\in {\cal T}}
V^*_{n-\lfloor n\tau \rfloor }(X)}{\inf_{\tau\in {\cal T}}V_{\lfloor n\tau\rfloor }(X)}.
\end{eqnarray}
This statistic has also the same form as statistic $R$ of \cite{Leybourne_Taylor_Kim:2007}, formed as a ratio of the minimized
CUSUMs of squared residuals obtained from the backward and forward subsamples of $X$, in the $I(0)/I(1)$ framework.
The limit distribution of these statistics is given in Proposition~\ref{consist}. To this end, define
\begin{equation}\label{Z*}
Z^*(u) \ := \ Z(1) - Z(1-u), \ u \in [0,1]
\end{equation}
and a continuous time analog of the partial sums' variance $V_{\lfloor n\tau\rfloor}(X)$ in \eqref{varS}:
\begin{eqnarray}
Q_\tau (Z)
&:=&\frac 1 {\tau^2}\bigg[\int_0^\tau \big(Z(u) - \frac{u}{\tau} Z(\tau)\big)^2 \d u -
\frac{1}{\tau}\Big(\int_0^\tau  \big(Z(u) - \frac{u}{\tau} Z(\tau)\big) \d u \Big)^2 \bigg]. \label{Q}
\end{eqnarray}
Note $Q_{1-\tau}(Z^*)$ is the corresponding analog of  $V^*_{n-\lfloor n\tau\rfloor}(X)$ in the numerators of the
statistics in (\ref{calL}) and ({\ref{RX}).

\begin{prop} \label{consist}
Assume ${\bf H_0}$. Then
\begin{eqnarray}\label{Dconv2}
\Big(\gamma^{-1}_{n} \big(S_{\lfloor n\tau_1\rfloor}  - \lfloor n\tau_1\rfloor A_{n}\big), \, \gamma^{-1}_{n}
\big(S^*_{\lfloor n\tau_2\rfloor} -
\lfloor n\tau_2\rfloor A_{n}\big)
\Big)&\longrightarrow_{D[0,1] \times D[0,
1]}&\big(Z(\tau_1), Z^*(\tau_2)\big).
\end{eqnarray}
Moreover, assume that
\begin{equation}\label{bddbelow}
Q_\tau (Z) >0 \quad \text {a.s.
for any} \quad \tau \in {\cal T}.
\end{equation}
Then
\begin{eqnarray}
\label{Wconv}
W_n(X)&\longrightarrow_{\rm law}&W(Z) \ := \ \sup_{\tau\in {\cal T}} \frac{Q_{1-\tau}(Z^*)}{Q_{\tau}(Z)},\nonumber   \\
I_n(X)&\longrightarrow_{\rm law}&I(Z) \ := \ \int_{\tau\in {\cal T}} \frac{Q_{1-\tau}(Z^*)}{Q_{\tau}(Z)}  \, \d \tau,   \\
R_n(X)&\longrightarrow_{\rm law}&R(Z) \ := \  \frac{\inf_{\tau\in {\cal T}} Q_{1-\tau}(Z^*)}{\inf_{\tau\in {\cal T}}Q_{\tau}(Z)}. \nonumber
\end{eqnarray}

\end{prop}

The convergence in \eqref{Dconv2} is an immediate consequence of ${\bf H_0}$, while the fact that \eqref{Dconv2} and \eqref{bddbelow} imply \eqref{Wconv} is a consequence of Proposition~\ref{general} stated in Section~\ref{sec:consistency}.

\begin{rem}
{\rm As noted previously, the alternative hypothesis ${\bf H_1}$ focuses on an
{\em increase} of $d$, and the statistics \eqref{memstat}, \eqref{RX} are defined
accordingly. It is straightforward to modify our testing procedures to test
for a decrease of persistence. In such case, the corresponding test statistics
are defined by exchanging forward and backward partial sums, or $V_{\lfloor n\tau\rfloor}(X)$
and $ V^*_{n-\lfloor n\tau\rfloor}(X)$:
\begin{eqnarray}\label{memstat*}\hskip-1cm
W^*_n(X)&:=&\sup_{\tau\in {\cal T}} {\cal L}^{-1}_n(\tau), \quad  \  I^*_n(X) \ := \ \int_{\tau\in {\cal T}} {\cal L}^{-1}_n(\tau) {\d} \tau, \quad
R^*_n(X)\ :=\ \frac{\inf_{\tau\in {\cal T}}
V_{\lfloor n\tau\rfloor}(X)}{\inf_{\tau\in {\cal T}}V^*_{n-\lfloor n\tau\rfloor}(X)}.
\end{eqnarray}
In the case when the direction of the change of $d$ is {\em unknown},
one can use various combinations of (\ref{memstat}) and (\ref{memstat*}), e.g. the sums
\begin{eqnarray*}
W^*_n(X) +W_n(X),   &I^*_n(X) +I_n(X),  &R^*_n (X) +R_n(X),
\end{eqnarray*}
or the maxima
\begin{eqnarray*}
\max\{W^*_n(X),W_n(X)\},   &\max\{I^*_n(X),I_n(X)\},  &\max\{R^*_n (X),R_n(X)\}.
\end{eqnarray*}
The limit distributions of the above six statistics under ${\bf H_0}$ follow
immediately from Proposition~\ref{consist}.
However, for a given direction of change, the `one-sided' tests in 
\eqref{memstat}, \eqref{RX} or \eqref{memstat*} are preferable as they are more powerful. }
\end{rem}

\subsection{Practical implementation for testing ${\bf H_0[I]}$  against ${\bf H_1[I]}$}
\label{sec:application-}

Under the `type I fBm null hypothesis' ${\bf H_0[I]}$, the limit distribution of the above statistics follows
from Proposition~\ref{consist} with $\gamma_n=n^{d+.5}$ and $Z= \kappa B^{\rm I}_{d + .5}$. In this case,
condition (\ref{bddbelow}) is verified and we obtain the following result.

\begin{cor} \label{consist0}
Assume ${\bf H_0[I]}$. Then
\begin{eqnarray}
\label{Wconv0}
W_n(X)&\longrightarrow_{\rm law}&W(B^{\rm I}_{d+.5}),  \quad
I_n(X)\ \longrightarrow_{\rm law}\ I(B^{\rm I}_{d+ .5}), \quad
R_n(X)\ \longrightarrow_{\rm law}\ R(B^{\rm I}_{d+.5}).
\end{eqnarray}
\end{cor}

The process $B^{\rm I}_{d + .5}$ in (\ref{Wconv0}) depends on {\it unknown} memory parameter $d$,
and so do the upper $\alpha-$quantiles  of the r.v.'s in the right-hand sides of (\ref{Wconv0})
\begin{equation}\label{quant}
q^{[\rm{I}]}_T(\alpha, d) := \inf \{x: \P(T(B^{\rm I}_{d+.5})\le x) \ge 1-\alpha)\},
\end{equation}
where $T =  W, I, R$.
Hence, applying the corresponding test, the unknown parameter $d$ in \eqref{quant}
is replaced by a consistent estimator $\hat d$.
\smallskip

\noindent \underline{Testing procedure.} \ 
Reject ${\bf H_0[I]}$, if
\begin{eqnarray}\label{tests}
W_n(X) \, > \, q^{[\rm{I}]}_W(\alpha, \hat d), \qquad I_n(X)\, >\, q^{[\rm{I}]}_I(\alpha, \hat d), \qquad R_n(X)\, >\, q^{[\rm{I}]}_R(\alpha, \hat d),
\end{eqnarray}
respectively, where $\hat d$ is a weakly consistent estimator of $d$:
\begin{equation} \label{dconsist}
\hat d \, \longrightarrow_p \, d, \qquad n \to \infty.
\end{equation}

\smallskip

The fact that the replacement of $d$ by $\hat d$ in (\ref{tests}) preserves  asymptotic significance level $\alpha$  is guaranteed
by the continuity of the quantile functions
provided by Proposition~\ref{quantcont} below.

\begin{prop} \label{quantcont} Let $d \in (-.5,1.5)$, $d\ne .5$,
$\alpha \in (0,1)$ and let $\hat d$ satisfy (\ref{dconsist}). Then
$$
q^{[\rm{I}]}_T(\alpha, \hat d)\, \longrightarrow_p \, q^{[\rm{I}]}_T(\alpha,d), \qquad \text{for} \quad  T\, = \, W,\,  I, \,  R.
$$
\end{prop}

We omit the proof of the  above proposition since it follows the same lines as in the paper Giraitis et al.\ (\citeyear{giraitis2}, Lemma 2.1)
devoted to tests of stationarity based on the V/S statistic.

Several estimators of $d$ can be used in (\ref{tests}). See the review paper \cite{bardet_lang} for a discussion
of some popular estimators. In our simulations we use the Non-Stationarity Extended Local Whittle Estimator (NELWE)
of \cite{abadir-FELW}, which applies to both stationary ($|d|<.5$) and
nonstationary ($d>.5$) cases.

\begin{rem} \label{rem3.2}
{\rm The above tests can be straightforwardly extended to  $d>1.5$, $d\ne 2.5,3.5,\dots$,
provided some modifications. Note that  
type I fBm for such values of $d$
is defined by iterating the integral in
\eqref{fbmI} (see e.g.\ \cite{davidsonJ}). On the other hand, although type I fBm can be defined for
$d=.5, 1.5, \dots $ as well, 
these values are excluded from
our discussion for the  following reasons. Firstly, in such case
the normalization $\gamma_n$ of partial sums process of $I(d)$
processes is different from $n^{d + .5} $ and contains an additional
logarithmic factor, see \cite{liu}. Secondly and more importantly, for $d=.5$  the limit
process $Z(\tau) = B^{\rm I}_1(\tau) = \tau B^{\rm I}_1(1)$ is a random line, in which case
the limit statistic $Q_\tau (Z) $ in (\ref{Q}) degenerates to zero, see also
Remark  \ref{remQdeg} below.
}
\end{rem}

\section{Consistency and asymptotic power}\label{sec:consistency}

It is natural to expect that under alternative hypotheses ${\bf H_1}$ or ${\bf H_1[I]}$, all three statistics
$W_n(X)$, $I_n(X)$, $R_n(X)$ tend to infinity in probability, provided the testing interval ${\cal T}$ and the
degeneracy interval $[0, \upsilon_1]$ of forward partial sums are embedded: ${\cal T} \subset [0, \upsilon_1]$.
This is true indeed, see Proposition~\ref{general} (iii) below, meaning that our tests are consistent.
Moreover, it is of interest to determine the rate at which these statistics grow under alternative, or the asymptotic power.
The following Proposition~\ref{general} provides the
theoretical background to study the consistency of the tests. It also
provides the limit distributions of the test statistics
under ${\bf H_0}$ since Proposition~\ref{consist} is an easy
corollary of Proposition~\ref{general} (ii).

\begin{prop} \label{general} (i)
Let there exist 
$0 \le \upsilon_0 <  \upsilon_1 \le 1$ and normalizations $\gamma_{ni} \to \infty$ and $A_{ni}, \, i=1,2$ such that
\begin{eqnarray} \label{Dconvpower}
\Big(\gamma^{-1}_{n1} \big(S_{\lfloor n\tau_1\rfloor}  - \lfloor n\tau_1\rfloor A_{n1}\big), \, \gamma^{-1}_{n2} \big(S^*_{\lfloor n\tau_2\rfloor} -
\lfloor n\tau_2\rfloor A_{n2}\big)
\Big)
&\longrightarrow_{D[0, \upsilon_1] \times D[0, 1-\upsilon_0]}&(Z_1(\tau_1), Z_2(\tau_2)\big),
\end{eqnarray}
where $(Z_1(\tau_1), Z_2(\tau_2)\big)$    
is a two-dimensional random process having a.s.\ continuous trajectories
on $[\upsilon_0, \upsilon_1] \times [1-\upsilon_1,1-\upsilon_0]$.
Then
\begin{eqnarray}\label{VVlim}
\big((n/\gamma_{n1}^2) V_{\lfloor n\tau_1\rfloor}(X), \,  (n/\gamma^2_{n2}) V^*_{n -
\lfloor n\tau_2\rfloor}(X) \big) &\longrightarrow_{D(0, \upsilon_1] \times
D[\upsilon_0, 1)}& \big(Q_{\tau_1} (Z_1),
Q_{1-\tau_2}(Z_2)\big).
\end{eqnarray}
Moreover, the limit process $\big(Q_{\tau_1} (Z_1),
Q_{1-\tau_2}(Z_2)\big)$ in (\ref{VVlim})
is a.s.\ continuous on $(\upsilon_0, \upsilon_1] \times
[\upsilon_0, \upsilon_1)$.

\medskip

\noi (ii)  Assume, in addition to (i), that ${\cal T} \subset {\cal U} := [\upsilon_0, \upsilon_1]$ and
\begin{equation}\label{bddbelow1}
Q_\tau (Z_1) >0 \quad \text {a.s.
for any} \quad \tau \in {\cal T}.
\end{equation}
Then, as $n \to \infty $,
\begin{eqnarray}
\label{WZZconv}
(\gamma_{n1}/\gamma_{n2})^2 W_n(X)&\longrightarrow_{\rm law}&\sup_{\tau\in {\cal T}} \frac{Q_{1-\tau}(Z_2)}{Q_{\tau}(Z_1)},\nonumber \\
(\gamma_{n1}/\gamma_{n2})^2 I_n(X)&\longrightarrow_{\rm law}&\int_{\tau\in {\cal T}} \frac{Q_{1-\tau}(Z_2)}{Q_{\tau}(Z_1)}  \, \d \tau,
\\
(\gamma_{n1}/\gamma_{n2})^2 R_n(X)&\longrightarrow_{\rm law}&\frac{\inf_{\tau\in {\cal T}} Q_{1-\tau}(Z_2)}{\inf_{\tau\in {\cal T}}Q_{\tau}(Z_1)}.
\nonumber \label{Rconv}
\end{eqnarray}

\medskip

\noi (iii) Assume, in addition to (i), that ${\cal T} \subset {\cal U}$,
 $Z_1(\tau) \equiv 0$, $\tau \in {\cal T}$
and the process $\{Q_{1-\tau} (Z_2), \tau \in {\cal T}\}$  is nondegenerate.
Then
\begin{eqnarray}\label{Winfty}
(\gamma_{n1}/\gamma_{n2})^2\left\{
\begin{array}{c}
W_n(X) \\
I_n(X) \\
R_n(X)
\end{array}
 \right\}&\longrightarrow_{p}&\infty.
\end{eqnarray}
\end{prop}

\medskip

\begin{rem} {\rm Typically, under ${\bf H_1}$ relation
 (\ref{Dconvpower}) is satisfied with $\gamma_{n2}$ increasing much faster than
 $\gamma_{n1}$ (e.g., $\gamma_{ni} = n^{d_i+.5}$, $i=1,2$, $d_1 < d_2$) and then (\ref{WZZconv})
imply that $W_n(X)$, $I_n(X)$ and $R_n(X)$ grow as  $O_{p}\big((\gamma_{n2}/\gamma_{n1})^2\big)$.
Two classes of fractionally integrated series with
changing memory parameter and satisfying  (\ref{Dconvpower}) are discussed in
Section~\ref{sec:appl-fract-integr}.
}
\end{rem}

\begin{rem} \label{remQdeg}
 {\rm Note that $Q_\tau (Z) \ge 0$ by the Cauchy-Schwarz inequality and that  $Q_\tau (Z) = 0$ implies
$Z(u) - \frac{u}{\tau} Z(\tau) = a $  \,  for all $u \in [0, \tau]$  and  some (random) $a = a(\tau)$.
In other words, $\P(Q_{\tau} (Z) = 0) >0$ implies that for some (possibly, random) constants
$a $ and $b$,
\begin{equation} \label{Z1deg}
\P \big(Z(u) = a + \frac{u}{\tau}\, b, \, \forall \, u \in [0, \tau] \big) \ > \ 0.
\end{equation}
Therefore, condition (\ref{bddbelow1}) implicitly excludes situations as in (\ref{Z1deg}), with $a\ne0$, $b\ne 0$,
which may arise under the null hypothesis ${\bf H_0}$, if  $A_n=0$ in (\ref{Dconv1}) whereas the $X_j$'s have
nonzero mean.
}
\end{rem}

\begin{rem}\label{kpss}
All the results in Sections~\ref{sec:testing-procedure}
and \ref{sec:appl-fract-integr} hold for Kim's statistics in (\ref{kim_integral}),
defined by replacing $V_{\lfloor n\tau\rfloor}(X), V^*_{n-\lfloor n\tau\rfloor}(X)$ in (\ref{memstat}),
(\ref{RX}) by $U_{\lfloor n\tau\rfloor}(X), U^*_{n-\lfloor n\tau\rfloor}(X)$ as given in (\ref{kim00}),
with the only difference that the functional $Q_\tau (Z)$  in the corresponding statements
must be replaced by its counterpart $\widetilde Q_\tau (Z) := \tau^{-2}\int_0^\tau \big(Z(u) -
\frac{u}{\tau} Z(\tau)\big)^2 \d u$, cf.\ (\ref{Q}).
\end{rem}

\section{Application to fractionally integrated processes}
\label{sec:appl-fract-integr}
This section discusses the convergence of forward and backward partial sums for some fractionally integrated
models with constant or changing memory parameter and the behavior of statistics $W_n, I_n, R_n$
for such models.

\subsection{Type I fractional Brownian motion and the null hypothesis ${\bf H_0[I]}$}
\label{sec:type-i-fractional}

It is well-known that type I fBm  arises
in the scaling limit of $d-$integrated, or $I(d)$, series with i.i.d.\ or martingale difference innovations.
See \cite{Davydov-70}, \cite{peligrad}, \cite{marinucci}, \cite{bruzaite} and the references
therein.

A formal definition of $I(d)$ process (denoted $\{X_t \} \sim I(d)$)
for $d > -.5$, $d \ne .5,1.5,\dots$ is given below. Let ${\rm MD}(0,1)$
be the class of all stationary ergodic martingale differences $\{\zeta_s, s \in \Z\}$
with unit variance $\E [\zeta_0^2] = 1$ and zero conditional expectation
$\E [\zeta_s |{\cal F}_{s-1}] = 0, \, s\in \Z$, where $\{{\cal F}_s, s \in \Z \}$ is a nondecreasing
family of $\sigma-$fields.

\begin{defn} \label{defId}
(i) Write  $\{X_t \} \sim I(0)$ if
\begin{equation}\label{I00}
 X_t = \sum_{j=0}^\infty a_j \zeta_{t-j}, \qquad t\in \Z
\end{equation}
 is a moving average with martingale difference innovations $\{\zeta_j \} \in {\rm MD}(0,1)$
 and summable coefficients $\sum_{j=0}^\infty |a_j | < \infty$, $\sum_{j=0}^\infty a_j \ne 0$.
\smallskip

\noindent (ii) Let $d \in (-.5,.5)\backslash\{0\}$. Write  $\{X_t\}\sim I(d)$ if $\{X_t\}$
is a fractionally integrated process
\begin{equation}\label{Xii}
 X_t = (1 -L)^{-d} Y_t = \sum_{j=0}^\infty \pi_j(d) Y_{t-j},  \qquad t \in \Z,
\end{equation}
 where $Y_t= \sum_{j=0}^\infty a_j \zeta_{t-j}$, $\{Y_t\} \sim I(0)$ and $\{\pi_j(d), j \ge 0\}$
 are the coefficients of the binomial expansion $(1-z)^{-d}=\sum_{j=0}^\infty\pi_j(d) z^j$, $|z|<1$.

\smallskip

\noindent (iii) Let $d >.5 $ and $d\ne 1.5,\, 2.5,\dots$. Write  $\{X_t\} \sim I(d)$ if
$X_t = \sum_{j=1}^t Y_j, \, t = 1,2, \dots$, where $\{Y_t\} \sim I(d-1)$.
\end{defn}

In the above definition, $\{X_t\} \sim I(d)$
for  $d>.5$ is recursively defined for $t=1,2,\dots$ only, as a $p-$times
integrated stationary $I(d-p)$ process, where $p = \lfloor d + .5 \rfloor $ is the integer part of $d + .5$,
and therefore $\{X_t\}$ has stationary increments of order $p$. A related definition of
$I(d)$ process involving initial values $X_{-i}, i =0,1,\dots $ is given in (\ref{Xiii})  below.
From Definition~\ref{defId}  it also follows that an $I(d)$ process can be written as a weighted sum
of martingale differences $\{\zeta_s \}\in {\rm MD}(0,1)$, for instance:
\begin{eqnarray} \label{Id}
X_t&=&\begin{cases}
\sum_{s \le t} (a \star \pi(d))_{t-s} \zeta_{s}, &-.5 < d < .5, \\
\sum_{s \le t} \sum_{1 \vee s\le j \le t} (a \star \pi(d-1))_{j-s} \zeta_{s}, &  .5<d<1.5,
\end{cases}\qquad t=1,2, \dots,
\end{eqnarray}
where
$(a \star \pi(d))_{j} := \sum_{i=0}^j a_i \pi_{j-i}(d)$, $j\ge 0$ is the convolution of
the sequences $\{a_j\}$ and $\{\pi_j(d)\}$.

\begin{prop} \label{linear} (i) Let $\{X_t\} \sim I(d)$ for some  $d\in (-.5,1.5)$, $d\ne .5$. If $d \in (-.5, 0]$, assume in addition
$\E |\zeta_1|^p < \infty, $ for some $p> 1/(.5 + d)$.
Then (\ref{Dconv00}) holds with $A_n = 0, \, \kappa =  \sum_{i=0}^\infty a_i$.
\smallskip

\noindent (ii) Let $\{\sigma_s, s \in \Z\} $ be an almost periodic sequence
such that $\bar \sigma^2 := \lim_{n\to \infty} n^{-1}\sum_{s=1}^n \sigma^2_s > 0$.
Let $\{X_t\}$ be defined as in (\ref{Id}), where $\, \zeta_s, \, s \in \Z \,$ are replaced by $ \, \sigma_s \zeta_s,\, s \in \Z \, $ and where $d$ and $\{\zeta_s\}$ satisfy the conditions in (i). Then (\ref{Dconv00}) holds with $A_n = 0, \kappa =  \bar \sigma \sum_{i=0}^\infty a_i. $
\end{prop}

The proof of Proposition~\ref{linear} can be easily reduced to the case $a_j=\kappa \delta_j$,
where $\delta_j = \1 (j=0)$ is the delta-function. Indeed,
\begin{equation} \label{approxF}
\E \bigg(\sum_{j=1}^n (X_j -X^\dagger_j) \bigg)^2 =  o\big(n^{2d+1}\big),
\end{equation}
where $X^\dagger_j:= \kappa (1- L)^{-d} \zeta_j $ \, $(-.5 < d < .5)$ and
$X^\dagger_j:=  \kappa \sum_{k=1}^j (1-L)^{-(d-1)} \zeta_k $ $(.5 < d < 1.5)$ is
(integrated) FARIMA$(0,d,0)$ process. The proof of the approximation (\ref{approxF})
is given in Section~\ref{sec:app}. The proof of Proposition~\ref{linear}
is omitted in view of (\ref{approxF}) and since similar results under slightly different
hypotheses on the innovations $\{\zeta_s\}$ can be found in \cite{bruzaite},
\cite{chan-terrin}, \cite{davidsonJ},  \cite{giraitis3}, and elsewhere.
In particular, the proof of the tightness in $D[0,1]$ given in \citet[Proposition 4.4.4]{giraitis3}
carries over to martingale difference innovations, see also \citet[Theorem 1.2]{bruzaite},
while part (ii) follows similarly to \citet[Theorem 1.1]{bruzaite}, using the fact that the sequence $\{\sigma_s \zeta_s \}$ satisfies
the martingale central limit theorem:
$n^{-1/2} \sum_{s=1}^{\lfloor n \tau\rfloor} \sigma_s \zeta_s\longrightarrow_{fdd} \bar\sigma B(\tau)$.
Note that the linear process $\{X_t\}$ in Proposition~\ref{linear} (ii) with
heteroscedastic noise $\{\sigma_s \zeta_s \}$ is nonstationary even if $|d|<.5$.

\subsection{Type II fractional Brownian motion and the null hypothesis ${\bf H_0[II]}$}\label{sec:typeII}

\begin{defn} \label{FII}  A type II fractional Brownian motion with parameter $d>-.5$ is defined by
\begin{eqnarray}\label{fbmII}
 B^{\rm II}_{d+.5}(\tau)&:=&\frac{1}{\Gamma(d+1)}\int_0^\tau (\tau -u)^d \d B(u),  \qquad \tau \ge 0,
\end{eqnarray}
where $\{B(u),u\ge 0\}$ is a standard Brownian motion with
zero mean and variance $\E B^2(u)=u$.
\end{defn}

A type II fBm shares many properties of type I fBm except that it has nonstationary increments, however, for
$|d| <.5 $ increments at time $\tau$ of type II fBm tend to those of type I fBm when $\tau\to \infty$.
\cite{davidsonH} discussed distinctions between the distributions of type I and type II fBms. Convergence
to type II fBm of partial sums of fractionally integrated processes was studied in
\cite{marinucci}. See also \cite{marinucci2}, \cite{davidsonJ},  \cite{ls}.

Type II fBm may serve as the limit process in the following specification of
the null hypothesis ${\bf H_0}$.

\medskip

\noindent ${\bf H_0[II]}$: {\it There exist $d>-.5$, $\kappa >0$ and a normalization $A_n$ such that
\begin{eqnarray}\label{DconvII}
n^{-d -.5} \big(S_{\lfloor n\tau\rfloor} - \lfloor n\tau\rfloor A_n \big)
&\longrightarrow_{D[0,1]}& \kappa B^{\rm II}_{d+.5} (\tau).
\end{eqnarray}
}

The alternative hypothesis to ${\bf H_0[II]}$ can be
again  ${\bf H_1[I]}$ 
of Section~\ref{sec:hyp}.

Proposition~\ref{linear} can be extended to type II fBm convergence in (\ref{DconvII}) as follows.
Introduce a `truncated' $I(0)$ process
\begin{equation}\label{YII}
Y_t := \begin{cases} \sum_{j=0}^t a_j \zeta_{t-j}, &t=1,2, \dots, \\
0, &t=0,-1,-2, \dots,
\end{cases}
\end{equation}
where $\{a_j\}$ and $\{\zeta_s\}$ are the same as in (\ref{I00}).
Following \cite{joh}, for $d > 0$ consider a $d-$integrated process  $\{X_t, t=1,2, \dots\}$ with given initial
values $\{X^0_{-i}, i=0,1,\dots \}$ as defined by
\begin{eqnarray}\label{Xiii}
 X_t&=&(1-L)_+^{-d} Y_t \ + \ (1-L)_+^{-d} (1-L)_-^d X^0_t, \qquad t=1,2, \dots,
\end{eqnarray}
where $\{Y_t\}$ is defined in (\ref{YII}) and the operators $(1-L)^d_{\pm}$
are defined through corresponding `truncated' binomial expansions:
$$
(1-L)^d_+ Z_t \ :=\ \sum_{j=0}^{t-1} \pi_j(-d) L^j Z_t, \qquad
(1-L)^d_- Z_t\ :=\ \sum_{j=t}^\infty \pi_j(-d) L^j Z_{t}\ =\ \sum_{i=0}^\infty\pi_{t+i}(-d) L^i Z_{0},
$$
$t=1,2,\dots$. Note that the term $(1-L)_+^{-d} (1-L)_-^d X^0_t$ in (\ref{Xiii}) depends
on initial values $\{X^0_{-i}, i = 0,1,\dots \}$ only. The choice of zero
initial values $X^0_{-i}=0$, $i=0,1,\dots$ in (\ref{Xiii}) leads to type
II process $X_t=(1-L)_+^{-d} Y_t$, \, more explicitly,
\begin{eqnarray} \label{IId}
X_t&=&
\sum_{s=1}^t  (a \star \pi(d))_{t-s} \zeta_{s}.
\end{eqnarray}
In general, $\{X^0_{-i}\}$ can be deterministic or random variables satisfying mild boundedness conditions
for the convergence of the series $(1-L)_-^d X^0_t$.

\begin{prop} \label{linearII}
(i) Let $\{X_t\}$ be defined in (\ref{Xiii}), with $\{Y_t\}$ as in (\ref{YII})
and initial values $\{X^0_{-i}\}$ satisfying for $d>.5$
\begin{equation}\label{initial}
\sup_{i \ge 0}\E (X^0_{-i})^2 <  \infty. 
\end{equation}
For $-.5 < d \le .5$ assume that $X^0_{-i} \equiv 0$. If $d \in (-.5, 0]$, assume in addition
$\E |\zeta_1|^p < \infty$, for some $p> 1/(.5 + d)$. Then (\ref{DconvII}) holds with $A_n = 0, \, \kappa = \sum_{i=0}^\infty a_i$.

\smallskip

\noindent (ii) Let $\{\sigma_s, s \ge  1\} $ be an almost periodic sequence
such that $\bar \sigma^2 := \lim_{n\to \infty} n^{-1}\sum_{s=1}^n \sigma^2_s > 0$.
Let $\{X_t\}$ be defined as in (\ref{IId}), where $\, \zeta_s, \, s \ge 1$ are replaced by $ \, \sigma_s \zeta_s,\, s \ge 1 \, $ and where $d$ and $\{\zeta_s\}$ satisfy
the conditions in (i). Then (\ref{DconvII}) holds with $A_n = 0$, $\kappa = \bar \sigma \sum_{i=0}^\infty a_i$.
\end{prop}

\begin{rem} For $d > .5$, Proposition~\ref{linearII} (i) implies that any
$L^2-$bounded initial values have no effect on the limit distribution of partial sums of the
process in (\ref{Xiii}). As it follows from the proof in Section~\ref{sec:app} below,
the above statement also remains valid for arbitrary initial values  $\{X^0_{-i} \}$ possibly
depending on $n$ and growing at a rate $O_p(n^{\lambda/2}) $ with some $0< \lambda < 1 \wedge (2d-1)$, viz.,
$\sup_{i \ge 0} \E (X^0_{-i})^2 <  C n^{\lambda}$, $d >0$.
\end{rem}

\smallskip

Similarly to Corollary~\ref{consist0}, Proposition~\ref{consist} implies the following
corollary.


\begin{cor} \label{corH0} Let  $\{X_t\}$ satisfy the conditions of Proposition~\ref{linearII}. Then
\begin{eqnarray}\label{WlimII}
W_n(X)&\longrightarrow_{\rm law}&W(B^{\rm II}_{d+.5}), \quad
I_n(X)\ \longrightarrow_{\rm law}\  I(B^{\rm II}_{d+.5}), \qquad
R_n(X)\  \longrightarrow_{\rm law}\  R(B^{\rm II}_{d+.5}),
\end{eqnarray}
where $\{B^{\rm II}_{d+.5}(\tau), \tau \in [0,1]\}$ is a type II fBm
as defined in (\ref{fbmII}).

\end{cor}


\begin{rem}
 {\rm Numerical experiments confirm that the upper
 quantiles $q^{[\rm{II}]}_T(\alpha,d)$, $T=W,I,R$, of the limit r.v.s on the r.h.s.\ of
 (\ref{WlimII}) are very close to the corresponding upper quantiles $q^{[\rm{I}]}_T(\alpha,d)$
 of the limiting statistics in (\ref{Wconv0}) when $d$ is smaller than 1 (see
 Figure~\ref{fig:I-II} in the particular case $T=I$). In other words, from a practical point of view, there
 is not much difference between type I fBM and type II fBm  null
 hypotheses ${\bf H_0[I]}$ and ${\bf H_0[II]}$ in testing for a change of $d$ when $d<1$.}
\end{rem}

\begin{figure}[h]
  \centering
\includegraphics[width=.7\textwidth,height=.3\textheight]{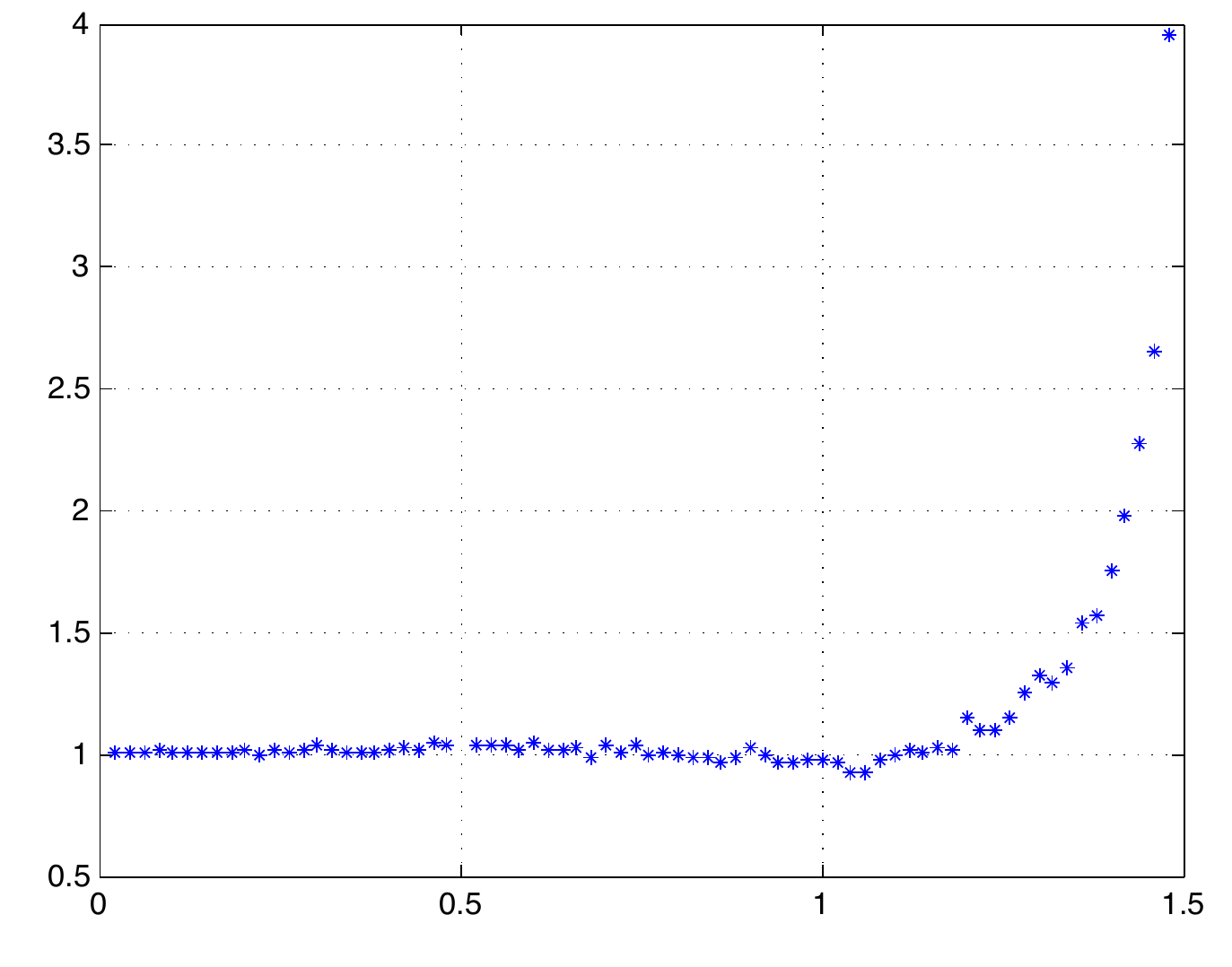}
  \caption{Representation of the ratio $q^{[\rm{I}]}_I(0.05,d) /
    q^{[\rm{II}]}_I(0.05,d)$ as function of $d$, with the choice $\underline\tau= 0.1$ and $\overline\tau=0.9$.}
\label{fig:I-II}
\end{figure}

\subsection{Fractionally integrated models with changing memory parameter}
\label{sec:fract-integr-models}



Let us discuss two nonparametric classes of nonstationary
time series with time-varying long memory parameter termed
`rapidly changing memory'  and `gradually changing memory'.
\smallskip

\noindent {\it Rapidly changing memory.} This class is obtained
by replacing parameter $d$ by a function $d(t/n)  \in [0,\infty)$ in the
FARIMA$(0,d,0)$ filter  
\begin{eqnarray} \label{farima}
\pi_j(d)&=&\frac{d}{1} \cdot \frac{d + 1}{2} \cdots
\frac{d - 1 + j}{j} \ = \
\frac{\Gamma(d + j)}{j! \Gamma (d)},\quad j=1,2,\dots, \quad \pi_0(d):=1.
\end{eqnarray}
Let $d(\tau), \, \tau \in [0,1]$ be a function taking values in the interval
$[0, \infty)$. (More precise conditions on the function  $d(\tau)$ will be specified below.)
 Define
\begin{eqnarray}
&&b_{1,j}(t)\ :=\  
\pi_{j}\big(d(\mbox{$\frac{t}{n}$})\big), \quad j=0,1,\dots,
\nonumber \\
&&X_{1,t}\ :=\  \sum_{s=1}^t b_{1,t-s}(t) \zeta_{s}, \qquad t=1, \dots, n, \label{X1t}
\end{eqnarray}
where the innovations $\zeta_s, s \ge 1$ satisfy the conditions of Definition~\ref{defId}. The
particular case
\begin{eqnarray}\label{d01}
d(\tau) = \begin{cases}
0, &\tau \in [0, \theta^*], \\
1, &\tau \in (\theta^*, 1]
\end{cases}
\end{eqnarray}
for some $0< \theta^* < 1$, leads to the model
\begin{eqnarray}
X_{1,t}&=&\begin{cases}
\zeta_t, &t = 1,2,\dots,  \lfloor \theta^* n\rfloor, \\
\sum_{s=1}^t \zeta_s, &t= \lfloor \theta^* n\rfloor +1, \dots, n,
\end{cases}\label{X10}
\end{eqnarray}
which corresponds to transition $I(0)\to I(1)$ at time $\lfloor \theta^* n\rfloor +1$.
A more general step function
\begin{eqnarray}\label{djump}
d(\tau) = \begin{cases}
d_1, &\tau \in [0, \theta^*], \\
d_2, &\tau \in (\theta^*, 1]
\end{cases}
\end{eqnarray}
corresponds to $\{X_{1,t}\} $ changing from $I(d_1)$ to $I(d_2)$ at time $\lfloor\theta^* n\rfloor+1$.

\smallskip

\noindent {\it Gradually changing memory.}
This class of nonstationary time-varying fractionally integrated processes
was defined in \cite{psv3,psv1,psv2}. Here, we use a truncated modification
of these processes with slowly varying memory parameter $d(t/n)\in [0,\infty)$, defined as
\begin{eqnarray}
 b_{2,j}(t) &:=&  \frac{d(\frac{t}{n})}{1} \cdot \frac{d(\frac{t-1}{n}) + 1}{2}  \cdots
 \frac{d(\frac{t-j+1}{n}) - 1 + j}{j}, \quad j = 1,2, \dots, \quad b_{2,0}(t):=1,\nonumber \\
 X_{2,t} & :=& \sum_{s=1}^t b_{2,t-s}(t) \zeta_s, \quad t=1, \dots, n.  \label{X2t}
\end{eqnarray}
Contrary to (\ref{X1t}), the process in (\ref{X2t}) satisfies an autoregressive time-varying fractionally integrated
equation with $\zeta_t$ on the right-hand side, see \cite{psv2}.
In the case when $d(\tau) \equiv d$ is constant function, the coefficients $b_{2,j}(t)$ in
(\ref{X2t}) coincide with FARIMA$(0,d,0)$ coefficients in (\ref{farima}) and in this case the
processes $\{X_{1,t}\}$ and $ \{X_{2,t}\}$ in (\ref{X1t}) and (\ref{X2t}) coincide. 

To see the difference between these two classes, consider the case of step function in  (\ref{d01}). Then
\begin{eqnarray}
X_{2,t}&=&\begin{cases}
\zeta_t, & t = 1,2, \dots, \lfloor \theta^* n\rfloor, \\
\sum_{s=\lfloor \theta^* n\rfloor+1}^t \zeta_s  + \sum_{s=1}^{\lfloor \theta^* n\rfloor}\frac{t - \lfloor \theta^*n\rfloor}{t-s} \zeta_s,
 &t= \lfloor \theta^* n\rfloor+1, \dots, n.
\end{cases} \label{X20}
\end{eqnarray}
Note $\frac{t - \lfloor \theta^*n\rfloor}{t-s} = 0 $ for  $t= \lfloor \theta^* n\rfloor$ and monotonically increases with $t \ge \lfloor \theta^* n\rfloor$. Therefore,
(\ref{X20}) embodies a gradual transition from $I(0)$ to $I(1)$, in contrast to an abrupt change of these regimes in (\ref{X10}). The
distinction between the two models (\ref{X10}) and (\ref{X20}) can be clearly seen from the variance behavior:
the variance of $X_{1,t}$ exhibits a jump from $1$ to $\lfloor \theta^* n\rfloor +1 = O(n)$ at time $t=\lfloor \theta^* n\rfloor +1 $,
after which it linearly increases with $t$,
while the variance of $X_{2,t}$ changes `smoothly' with $t$:
\begin{eqnarray*}
{\rm Var}(X_{2,t})&=&\begin{cases}
1, &t = 1,2, \dots,\lfloor \theta^* n\rfloor, \\
(t - \lfloor \theta^*n\rfloor)\, + \,\sum_{s=1}^{\lfloor \theta^* n\rfloor}\frac{(t - \lfloor \theta^*n\rfloor)^2}{(t-s)^2}, &t= \lfloor \theta^* n\rfloor +1, \dots, n.
\end{cases}
\end{eqnarray*}
Similar distinctions between  (\ref{X1t}) and (\ref{X2t}) prevail also in the case of general `memory function'  $d(\cdot)$: when the memory parameter
$d(t/n)$ changes
with $t$, this change gradually affects the lagged ratios in the coefficients $b_{2,j}(t)$ in (\ref{X2t}), and not all lagged
ratios simultaneously as in the case of $b_{1,j}(t)$, see (\ref{farima}).

\subsection{Asymptotics of change-point statistics for fractionally integrated models with
changing memory parameter}\label{sec:5.4}

In this subsection we study the joint convergence of forward and backward partial sums as in (\ref{Dconv}) for
the two models in (\ref{X1t}) and  (\ref{X2t}) with time-varying memory parameter $d(t/n)$. After the statement
of Proposition~\ref{power} below, we discuss its implications for the asymptotic power of our tests.

\vskip.1cm

Let us specify a class of `memory functions' $d(\cdot)$.  For $0< d_1 < d_2 < \infty$ and $0\leq \underline\theta \le \overline\theta \leq 1$, introduce the class ${\cal D}_{\underline\theta,\overline\theta}(d_1, d_2)$ of left-continuous nondecreasing  functions $d (\cdot) \equiv \{d (\tau), \tau \in [0,1]\}$ such that
\begin{eqnarray}\label{dchange}
d(\tau) = \begin{cases}
d_1, &\tau \in [0, \underline\theta], \\
d_2, &\tau \in [\overline\theta, 1],
\end{cases},   \qquad d_1  < d(\tau) < d_2,  \quad  \underline\theta < \tau <\overline \theta.
\end{eqnarray}
The interval $\Theta := [\underline\theta,\overline\theta]$ will be called the {\it memory change interval}. 
Note that for $\underline\theta=\overline\theta \equiv \theta^*$, the class ${\cal D}_{\theta^*,\theta^*}(d_1, d_2)$ consists
of a single step function in (\ref{djump}). Recall from Section~\ref{sec:testing-procedure} that
the interval ${\cal T} = [\underline\tau,\overline\tau]$ in memory change statistics in (\ref{memstat}) and \eqref{RX} is called the (memory) testing interval.
When discussing
the behavior of memory tests under alternatives in (\ref{X1t}), (\ref{X2t}) with changing memory parameter,
the intervals $\Theta$ and ${\cal T}$ need not coincide since $\Theta$ is not known {\it a priori}.


With a given $d (\cdot) \in {\cal D}_{\underline\theta,\overline\theta}(d_1, d_2), $ 
we associate a 
function
\begin{equation} \label{Huv}
H(u,v):=
\begin{cases}
  \int_u^{v} \frac{d(x) - d_2}{v-x}\, \d x, & 0\le u \le v \le
  1,  \\ 0, & \text{otherwise}
\end{cases}
\end{equation}
Note $H(u,v) \le 0$ since $d(x) \le d_2, \, x \in [0,1]$ and $H(u,v) = 0$ if  $\overline\theta \le u \le v \le 1$.
Define two Gaussian processes ${\cal Z}_1 $ and ${\cal Z}_2 $ by
\begin{eqnarray} \label{Gd}
{\cal Z}_1(\tau)
&:=&\frac{1}{\Gamma(d_2)} \int_0^{\tau}  \Big\{\int_{\overline\theta}^\tau (v-u)_+^{d_2-1}  \d v \Big\} \d B(u) \ = \
B_{d_2+.5}^{\rm II}(\tau) -  B_{d_2+.5}^{\rm II}(\overline\theta),\nonumber \\
{\cal Z}_2(\tau)
&:=&\frac{1}{\Gamma(d_2)} \int_0^{\tau}  \Big\{\int_{\overline\theta}^\tau (v-u)_+^{d_2-1} \e^{H(u,v)} \d v \Big\} \d B(u), \quad \tau >\overline \theta,\\
{\cal Z}_1(\tau)
&=&{\cal Z}_2(\tau) \ := \ 0, \quad \tau \in [0, \overline\theta]. \nonumber
\end{eqnarray}
The processes $\{{\cal Z}_i(\tau), \in [0,1]\}$, $i=1,2$
are well-defined for any $d_2 > -0.5$ and have 
a.s.\ continuous trajectories.
In the case  $\underline\theta = \overline\theta \equiv \theta^*$ and a step function $d(\cdot)$ in (\ref{dchange}),
${\cal Z}_2(\tau)$ for $\tau > \theta^*$  can be rewritten as
\begin{eqnarray}\label{Z2cp}
{\cal Z}_2(\tau) &=&
\frac{1}{\Gamma(d_2)} \int_0^{\tau}  \Big\{\int_{\theta^*}^\tau (v-u)_+^{d_1-1} (v-\theta^*)^{d_2 - d_1} \d v \Big\} \d B(u).
\end{eqnarray}
Related class of Gaussian processes was discussed in \cite{psv2} and \cite{surg2}.

\begin{prop} \label{power} Let
$d (\cdot) \in {\cal D}_{\underline\theta, \overline\theta}(d_1, d_2)$ for some $0 \le d_1 < d_2 < \infty$,
$0 \leq \underline\theta \le\overline \theta \leq 1$. Let $S_{i,k}$ and $S^*_{i,n-k}$, $i=1,2$ be the
forward and backward partial sums processes corresponding to time-varying fractional filters
$\{X_{i,t}\}$, $i=1,2$ in (\ref{X1t}), (\ref{X2t}), with memory parameter $d(t/n)$ and standardized i.i.d.\ innovations $\{\zeta_j, \, j\ge 1\}$.  Moreover, in the case $d_1 = 0$ we assume that $\E |\zeta_1|^{2+\delta} < \infty$
for some $\delta >0$. Then

\noindent (i) 
for any $\theta\in (0,\underline\theta]$ with $\underline\theta>0$
\begin{eqnarray} \label{ZZ1}
\big(n^{-  d_1-.5} S_{i,\lfloor n\tau_1\rfloor}, n^{-d_2-.5} S^*_{i,\lfloor n\tau_2\rfloor}\big)
&\longrightarrow_{D[0, \theta] \times D[0, 1-\underline\tau]}&(Z_{i,1}(\tau_1), Z_{i,2}(\tau_2)\big),  \qquad  i=1,2,
\end{eqnarray}
where
\begin{eqnarray} \label{ZZij}
Z_{i,1}(\tau)&:=&B^{\rm II}_{d_1 + .5}(\tau), \qquad
Z_{i,2}(\tau)\, :=\, {\cal Z}^*_i(\tau) \, = \, {\cal Z}_i(1) - {\cal Z}_i(1-\tau), \qquad i=1,2,
\end{eqnarray}
and ${\cal Z}_i$, $i=1,2$ are defined in (\ref{Gd});
\smallskip

\noindent (ii) for any $\theta\in [\underline\theta,1]$,
for any $d > d(\theta), \,  d_1 < d < d_2$
\begin{eqnarray} \label{ZZ2}
\big(n^{-  d-.5} S_{i,\lfloor n\tau_1\rfloor}, n^{-d_2-.5} S^*_{i,\lfloor n\tau_2\rfloor}\big)
&\longrightarrow_{D\lfloor 0, \theta\rfloor \times D[0, 1-\underline\tau]}&(0, Z_{i,2}(\tau_2)\big),  \qquad  i=1,2,
\end{eqnarray}
where $Z_{i,2}, \, i=1,2$ are the same as in (\ref{ZZij}).
\end{prop}


The power of our tests depends on whether the testing and the memory change intervals  have
an empty intersection or not. When $\overline\tau<\underline\theta$, Proposition~\ref{power} (i) applies taking $\theta=\overline\tau$ and the asymptotic distribution of the memory test statistics for models (\ref{X1t}) and (\ref{X2t}) follows from
Proposition~\ref{general} (\ref{WZZconv}), with normalization $(\gamma_{n2}/\gamma_{n1})^2 = n^{2(d_1-d_2)} \to 0$, implying the consistency of the tests. But this situation is untypical for practical applications and hence not very interesting. Even less interesting seems the case  when a change of memory ends before the start of the testing interval, i.e., when $\overline\theta \le \underline\tau$. Although the last case is not covered by Proposition~\ref{power}, the limit distribution of the test statistics for models (\ref{X1t}), (\ref{X2t})  exists with trivial normalization
$(\gamma_{n2}/\gamma_{n1})^2 = 1$ and therefore our tests are inconsistent, which is quite natural in this case.

Let us turn to some more interesting situations, corresponding to the case when the intervals ${\cal T}$
and $\Theta$ have a nonempty intersection of positive length. There are two possibilities:


\smallskip

\noi {Case 1}: $\underline\tau < \underline\theta \leq \overline\tau$ (a change of memory occurs after the beginning of the testing interval), and
\smallskip

\noi {Case 2}: $\underline\theta\leq \underline \tau  < \overline\theta$ (a change of memory occurs before the beginning of the testing interval).
\smallskip

Let us consider Cases 1 and 2 in more detail.
\smallskip

\noi {\it Case 1.}  Let $\widetilde {\cal T} := [\underline\tau, \underline\theta] \subset {\cal T}$.
Introduce the following `dominated' (see \eqref{statdom}) statistics:
\begin{eqnarray}\label{widetildememstat}
\widetilde W_n(X)&:=&\sup_{\tau\in \widetilde {\cal T}} \frac{V^*_{n-\lfloor n\tau\rfloor}(X)}{V_{\lfloor n\tau\rfloor}(X)}, \quad  \
\widetilde I_n(X)\ := \ \int_{\widetilde {\cal T}} \frac{V^*_{n-\lfloor n\tau\rfloor}(X)}{V_{\lfloor n\tau\rfloor}(X)}\, \d \tau, \\
\widetilde R_n(X)&:=&\frac{\inf_{\tau\in {\cal T}} V^*_{n-\lfloor n\tau\rfloor}(X)}{\inf_{\tau\in \widetilde {\cal T}}V_{\lfloor n\tau\rfloor}(X)}.  \nonumber
\end{eqnarray}
Clearly,
\begin{equation}\label{statdom}
W_n(X) \ge \widetilde W_n(X), \qquad I_n(X) \ge \widetilde I_n(X), \qquad R_n(X) \ge \widetilde R_n(X), \qquad {\rm a.s.}
\end{equation}
The limit distribution of (\ref{widetildememstat}) for models (\ref{X1t}) and (\ref{X2t}) can be derived from
propositions~\ref{general} and \ref{power} (i) choosing $\theta=\underline\theta$. In particular, it follows that
$n^{2(d_1-d_2)}\widetilde W_n(X_i),\, n^{2(d_1-d_2)}\widetilde I_n(X_i)$,
and $n^{2(d_1-d_2)}\widetilde R_n(X_i), \, i=1,2$ tend, in distribution, to the corresponding limits in (\ref{WZZconv}),
with ${\cal T}$ replaced by $\widetilde {\cal T}$ and  $Z_1 = Z_{i,1},\, Z_2 = Z_{i,2}$, $i=1,2 $ as defined in
(\ref{ZZij}). Moreover, it can be shown that $n^{-2d_1} V_{\lfloor n \tau\rfloor}(X_i) \longrightarrow_p \infty $ for any
$\tau \in {\cal T}\,\backslash\,\widetilde {\cal T}$. Therefore, in Case 1, the limit distributions of the original statistics in
(\ref{memstat}) and the `dominated' statistics in (\ref{widetildememstat}) coincide.

\smallskip

\noi {\it Case 2.} In this case, define $\widetilde {\cal T} := [\underline\tau, \tilde \theta] \subset {\cal T}$, where
$\tilde \theta \in (\underline\tau,\overline \theta)$ is an inner point of the interval  $[\underline\tau,\overline \theta]$.
Let $\widetilde W_n(X)$, $\widetilde I_n(X)$, $\widetilde R_n(X)$ be defined as in (\ref{widetildememstat}).
Obviously, relations (\ref{statdom}) hold as in the previous case.
Since the memory parameter increases on the interval $\widetilde {\cal T}$, the limit distribution of the process
$V_{\lfloor n \tau\rfloor}(X_i)$, $\tau\in \widetilde {\cal T}$ in the denominator of the statistics is not
identified from Proposition \ref{general} (ii).
Nevertheless in this case we can use Propositions \ref{general} (iii) and \ref{power} (ii)
to obtain a robust rate of growth of the memory statistics
in (\ref{widetildememstat}) and (\ref{memstat}). Indeed from Proposition~\ref{power} (ii) with $\theta=\tilde\theta$,  we have that
$n^{-2d} V_{\lfloor n \tau\rfloor}(X_i) \longrightarrow_{D(0, \tilde \theta]} 0 $ for any $d_2 > d > d(\tilde \theta)$ and hence
$n^{2(d-d_2)}W_n(X_i)$, $n^{2(d-d_2)}I_n(X_i)$ and $n^{2(d-d_2)}R_n(X_i)$, $i=1,2$ tend to infinity, in probability.

\section{Testing in the presence of linear trend}\label{sec:lin-trend}

The tests discussed in Section~\ref{sec:testing-procedure} can be further developed to include
the presence of a linear trend. In such a case, partial sums $S_{\lfloor n\tau\rfloor}$ may
grow as a second-order polynomial of $\lfloor n\tau\rfloor$ (see Example~\ref{ex2} below).
Then, the null and alternative hypotheses have to be modified, as follows.
\medskip

\noindent ${\bf H}^{\rm trend}_0$: There exists normalizations $\gamma_{n} \to \infty, \, A_{n}, \, B_{n},$ such that
\begin{eqnarray}\label{Dtrendconv1}
 \gamma^{-1}_{n} \big(S_{\lfloor n\tau\rfloor} -  \lfloor n\tau\rfloor A_{n} - \lfloor n\tau\rfloor^2 B_{n} \big)
 &\longrightarrow_{D[0,1]}& Z(\tau),
\end{eqnarray}
where $\{Z(\tau), \tau \in [0,1]\}$ is a nondegenerate a.s.\
continuous random process.

\bigskip

\noindent ${\bf H^{\rm trend}_1}$: There exist $0\le\upsilon_0 <\upsilon_1 \le 1$ and normalizations
$\gamma_{n} \to \infty$, $A_n$, $B_n$, such that
\begin{eqnarray}\label{Dtrendconv}
&&\Big(\gamma^{-1}_{n} \big(S_{\lfloor n\tau_1\rfloor}-\lfloor n\tau_1\rfloor \,
 A_n - \lfloor n\tau\rfloor^2 B_{n} \big), \, \gamma^{-1}_{n}
 \big(S^*_{\lfloor n\tau_2\rfloor} -
 \lfloor n\tau_2\rfloor^* A_n - \lfloor n\tau_2\rfloor^{2*} B_n\big)\Big)\\
&&\longrightarrow_{D[0, \upsilon_1] \times D[0,1-\upsilon_0]}\ \big(0, Z_2(\tau_2)\big), \nonumber
\end{eqnarray}
where $\{Z_2(\tau), \tau \in [1-\upsilon_1,1-\upsilon_0]\}$ is a
nondegenerate  a.s.\ continuous random process, $\lfloor n\tau\rfloor^*:=n-\lfloor n(1 -\tau)\rfloor = \lfloor n\tau\rfloor, \,
\lfloor n\tau\rfloor^{2*} := n^2 - \lfloor n(1-\tau)\rfloor^2 $.

\smallskip

\begin{ex}\label{ex2}
Consider a process $\{X_t\}$ defined as in Example~\ref{ex1} from equation \eqref{changed}.
We construct the process $\{\mathcal{X}_t\}$ by adding to $\{X_t\}$ an additive linear trend:
\begin{equation}
 \mathcal{X}_t = X_t + a  + b t,\label{eq:Xtrend}
\end{equation}
where $a,b$ are some coefficients.

When $d_1=d_2=d$, we have $\mathcal{X}_t = \varepsilon_t(d)+a+bt$  and
$\{\mathcal{X}_t\}$ satisfies the hypothesis ${\bf H^{\rm trend}_0}$
with $B_{n} = \frac{b}{2}, A_{n} = a + \frac{b}{2}$ and $Z = B^{\rm I}_{d+ .5}$. Indeed,
\begin{eqnarray*}
&&n^{-d -.5} \Big(S_{\lfloor n\tau\rfloor}(\mathcal{X}) - \lfloor n\tau\rfloor (a + \frac{b}{2})  -  \lfloor n\tau\rfloor^2  \frac{b}{2} \Big)\\
&&=\  n^{-d -.5}S_{\lfloor n\tau\rfloor}(\veps(d)) +    n^{-d -.5}
\Big\{ a \lfloor n \tau\rfloor + b \sum_{j=1}^{\lfloor n\tau\rfloor} j
- \lfloor n\tau\rfloor \Big(a + \frac{b}{2}\Big) -\lfloor n\tau\rfloor^2 \frac{b}{2} \Big\} \\
&&= \ n^{-d -.5}S_{\lfloor n\tau\rfloor}(\veps(d)) \
\longrightarrow_{D[0,1]}\  B^{\rm I}_{d+.5} (\tau),\quad n\to\infty.
\end{eqnarray*}
\end{ex}

\smallskip

Under linear trend, the test statistics of Section~\ref{sec:test-statistics} have
to be modified, as follows. For a fixed $1 \le k \le n$, let
$$
\widehat X_j\ :=\  X_j - \hat a_k(X) - \hat b_k(X) j, \qquad 1 \le j \le k
$$
denote the residuals from the least-squares  regression of
$(X_j)_{1\le j \le k} $ on $(a + b j)_{1\le j \le k}$.
Similarly, let
$$
\widehat X^*_j\ := \ X_j - \hat a^*_{n-k}(X) - \hat b^*_{n-k}(X) j, \qquad k < j \le n
$$
denote the residuals from the least-squares regression of
$(X_j)_{k < j \le n}$ on $(a + b j)_{k < j \le n}$. The corresponding
intercept and slope coefficients are defined through
$(\hat a_k(X), \hat b_k(X))\ := \ {\rm argmin} \big(\sum_{j=1}^k (X_j - a - bj)^2 \big)$,
$(\hat a^*_{n-k}(X), \hat b^*_{n-k}(X)) \ := \ {\rm argmin} \big(\sum_{j=k+1}^n (X_j - a - bj)^2 \big)$.
The variance estimates of de-trended forward and backward partial sums are defined by
\begin{equation}\label{varStrend}
\begin{split}
{\mathcal V}_k(X)&\ :=\ \frac{1}{k^2} \sum_{j=1}^k (\widehat S_j)^2  - \bigg(\frac{1}{k^{3/2}} \sum_{j=1}^k \widehat S_j\bigg)^2,\\
{\mathcal V}^*_{n-k}(X)&\ :=\ \frac{1}{(n-k)^2} \sum_{j=k+1}^n \big(\widehat S^*_{n-j+1}\big)^2 -
\bigg(\frac{1}{(n-k)^{3/2}} \sum_{j=k+1}^n \widehat S^*_{n-j+1} \bigg)^2,
\end{split}
\end{equation}
where
$$
\widehat S_j \ :=\  \sum_{i=1}^j \widehat X_i = \sum_{i=1}^j (X_i - \hat a_k(X) - \hat b_k(X) i),
\quad  \widehat S^*_{n-j+1} \ :=\  \sum_{i=j}^n  \widehat X^*_i = \sum_{i=j}^n (X_i - \hat a^*_{n-k}(X) - \hat b^*_{n-k}(X) i),
$$
cf.\ (\ref{varS}).
Replacing $V_n(X), V^*_{n-k}(X)$ in (\ref{memstat})--(\ref{RX}) by the corresponding
quantities ${\mathcal V}_n(X), {\mathcal V}^*_{n-k}(X)$,
the statistics in presence of a linear trend are given by
\begin{equation}\label{trstat1}
\mathcal{W}_n(X)\ :=\  \sup_{\tau\in {\cal T}} \frac{{\mathcal
  V}^*_{n-\lfloor n\tau \rfloor }(X)}{{\mathcal
  V}_{\lfloor n\tau\rfloor }(X)}, \quad  \mathcal{I}_n(X) \ := \ \int_{\tau\in {\cal T}}\frac{{\mathcal
  V}^*_{n-\lfloor n\tau \rfloor }(X)}{{\mathcal
  V}_{\lfloor n\tau\rfloor }(X)}{\d} \tau,\quad  \mathcal{R}_n(X)\ :=\ \frac{\inf_{\tau\in {\cal T}}
{\mathcal V}^*_{n-\lfloor n\tau \rfloor }(X)}{\inf_{\tau\in {\cal T}}{\mathcal V}_{\lfloor n\tau\rfloor }(X)}.
\end{equation}
Note that (\ref{trstat1}) agree with (\ref{memstat})--(\ref{RX}) if no trend is
assumed (i.e. $b$ is known and equal to zero).

Under the null hypothesis ${\bf H^{\rm trend}_0}$ the distributions of (\ref{trstat1})
can be obtained similarly to that of (\ref{memstat})--(\ref{RX}). The following proposition
is the analog of the corresponding result (\ref{VVlim}) of  Proposition~\ref{general}~(ii)
for (\ref{varStrend}).

\begin{prop} \label{VVtrend} Under the hypothesis ${\bf H}^{\rm trend}_0$,
\begin{equation} \label{V1trend}
\big((n/\gamma_n^2) {\mathcal V}_{\lfloor n\tau_1\rfloor}(X), (n/\gamma_n^2) {\mathcal V}^*_{n-\lfloor n\tau_2\rfloor}(X)\big)
  \longrightarrow_{D(0, 1] \times
D[0, 1)} \big({\mathcal Q}_{\tau_1} (Z), {\mathcal Q}_{1-\tau_2}(Z^*)\big),
\end{equation}
where $Z^*$ is defined in (\ref{Z*}) and
\begin{eqnarray}\label{Qtrend}
{\mathcal Q}_\tau (Z)
&:=&\frac{1}{\tau^2} \bigg[\int_0^\tau {\mathcal Z}(u,\tau)^2 \d u - \frac{1}{\tau}
\Big(\int_0^\tau {\mathcal Z}(u,\tau)\d u\Big)^2\bigg], \\
{\mathcal Z}(u,\tau)&:=&Z(u) + 2Z(\tau)\frac{u}{\tau}\Big(1 -\frac{3u}{2\tau}\Big) +
6\Big(\frac{1}{\tau}\int_0^\tau Z(v) \d v \Big) \frac{u}{\tau} \Big(1-\frac{u}{\tau}\Big).  \nonumber
\end{eqnarray}

\end{prop}

Note that the process $\{{\mathcal Z}(u,\tau), u \in [0,\tau] \}$ defined in (\ref{Qtrend})
satisfies ${\mathcal Z}(0,\tau) = {\mathcal Z}(\tau,\tau) = 0$ and
$\int_0^\tau {\mathcal Z}(u,\tau)\d u = 0$. In the case of Brownian motion $Z = B$, $\{{\cal Z}(u,1),u\in [0,1]\}$
 is known as the second level
Brownian bridge (see \cite{macneill:1978}).
Extension of Proposition~\ref{consist} to the modified statistics $W_n, I_n, R_n$ with the ratio (\ref{calL}) replaced by
${\mathcal V}^*_{n-\lfloor n\tau\rfloor}(X)/{\mathcal V}_{\lfloor n\tau\rfloor}(X)$ is straightforward. Clearly,
${\cal Z}(u,\tau)$ in (\ref{Qtrend}) is different from the corresponding process
$Z(u) - \frac{u}{\tau} Z(\tau)$ in (\ref{Q}) and therefore the `de-trended' tests ${\mathcal W}_n, {\mathcal I}_n, {\mathcal R}_n$ have
different critical regions from those in (\ref{tests}). We also note that $\{{\cal Z}(u,\tau)\}$ can be heuristically
defined as the residual process $\{{\cal Z}(u,\tau) = Z(u) - \hat a_\tau  u - \frac{\hat b_\tau}{2} u^2\} $ from the
least squares regression of $(\d Z(u)/\d u)_{u \in [0,\tau]} $ onto $(a + b u)_{u \in [0,\tau]},$ with $(\hat a_\tau, \hat b_\tau)$
minimizing the integral $ \int_0^\tau (  \frac{\d Z(u)}{\d u} -   a - b u)^2 \d u $. Indeed, the above minimization
problem leads to linear equations $
\int_0^\tau (  \frac{\d Z(u)}{\d u} -   a - b u)  \d u = 0, \
\int_0^\tau (  \frac{\d Z(u)}{\d u} -   a - b u)u \d u = 0, $
or
$$
Z(\tau) =   a\tau + b\frac{\tau^2}{2}, \qquad
\tau Z(\tau) - \int_0^\tau Z(u) \d u\  = \   a \frac{\tau^2}{2} +  b \frac{\tau^3}{3},
$$
where we used $\int_0^\tau u \frac{\d Z(u)} {\d u} \d u = \int_0^\tau u \d Z(u) = \tau Z(\tau) - \int_0^\tau Z(u) \d u$. Solving
the above equations leads to the same $\hat a_\tau  = -\frac{2}{\tau} Z(\tau) + \frac{6}{\tau^2}
\int_0^\tau Z(u) \d u, \
\hat b_\tau = \frac{6}{\tau^2} Z(\tau) - \frac{12}{\tau^3}
\int_0^\tau Z(u) \d u $
as in (\ref{ablim}) below. The resulting expression of
${\cal Z}(u,\tau) = Z(u) - \hat a_\tau u - \frac{\hat b_\tau}{2} u^2 $ agrees with
(\ref{Qtrend}).

\section{Simulation study}\label{sec:simulation}

In this section we compare from numerical experiments
 the finite-sample performance of the three test statistics in \eqref{tests} for  testing
 $\bf{H_0[I]}$ against $\bf{H_1}[I]$ with nominal level $\alpha= 5\%$. A comparison with the
 Kim's tests based on the ratio \eqref{kim0}, is also provided.

The main
steps to implement the testing procedures defined in
\eqref{tests} are the following: 
\begin{itemize}
\item We choose $\overline\tau=1-\underline\tau$ for $\underline\tau\in(0,1)$
which  defines the testing region $\mathcal{T}:=[\underline\tau, 1-
  \underline\tau]$. Sensitivity to the choice of $\underline\tau$ is also explored;

\item For each simulated sample  $X_1,\dots,X_n,$ we estimate the parameter $d$
using the
  NELWE of \cite{abadir-FELW} as the estimate of $d$.
  Following the recommendation in \cite{abadir-FELW},
  the bandwidth parameter in the above  estimate is chosen  to be  $\lfloor \sqrt{n}\rfloor$;

\item 
The quantiles
$q_T^{[\rm{I}]}(.05, d)$ in the critical regions \eqref{tests}, as functions of $d$, for $T=W, R, I,$  and for chosen values
of $\underline\tau$, are approximated by extensive Monte Carlo experiments. The integral appearing in the definition
of $T=I$ in \eqref{Wconv} is approximated by a Riemann sum. See also \cite{hassler_scheithauer:2008} on approximation of similar
quantiles. The quantile graph for $T=I$ and three different values of $\underline\tau$ is shown in
Figure~\ref{fig:qq}.
\end{itemize}

\begin{figure}[h]
 \centering
\includegraphics[width=.57\textwidth,height=.25\textheight]{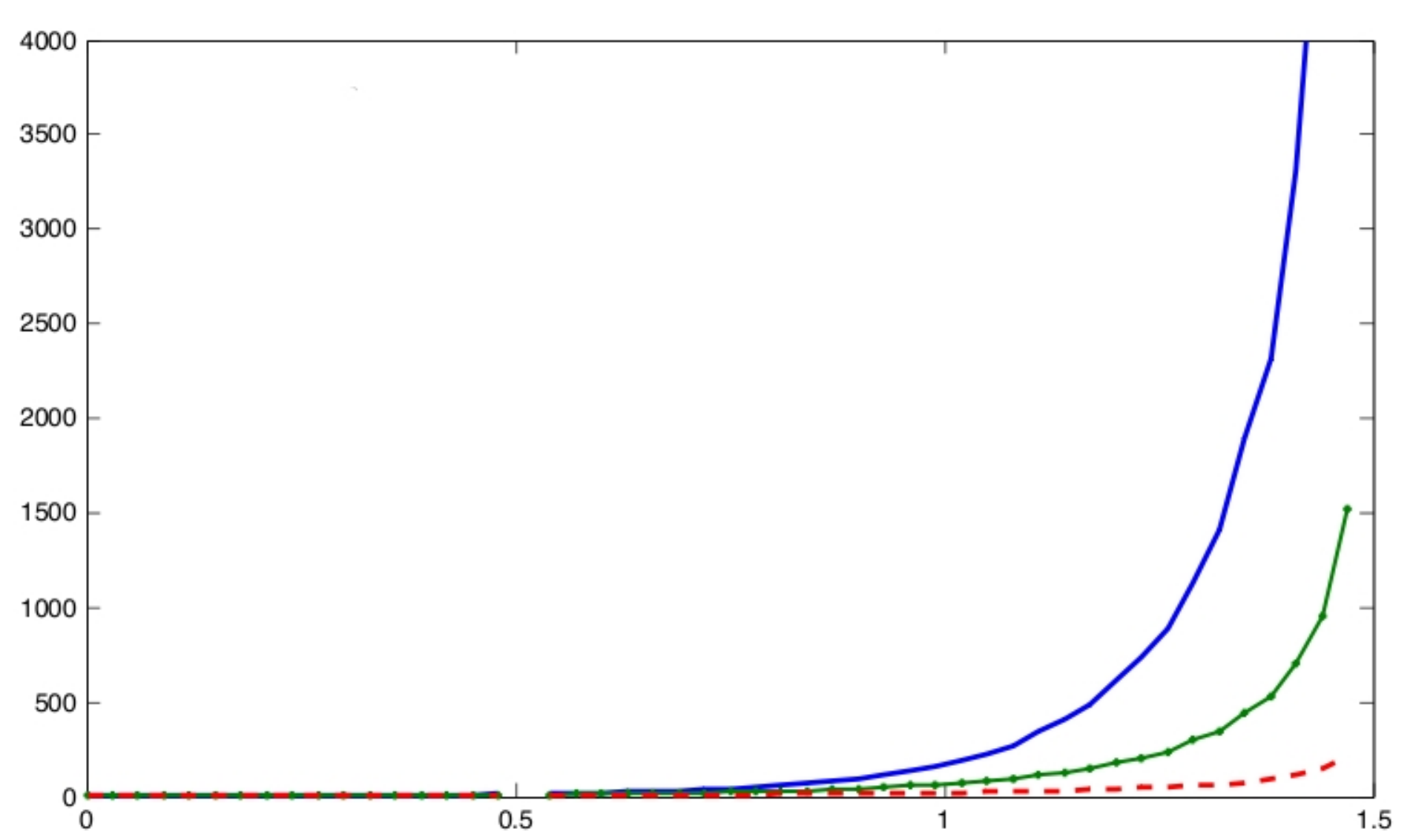}
\raisebox{-0.2cm}
{\includegraphics[width=.4\textwidth,height=.25\textheight]{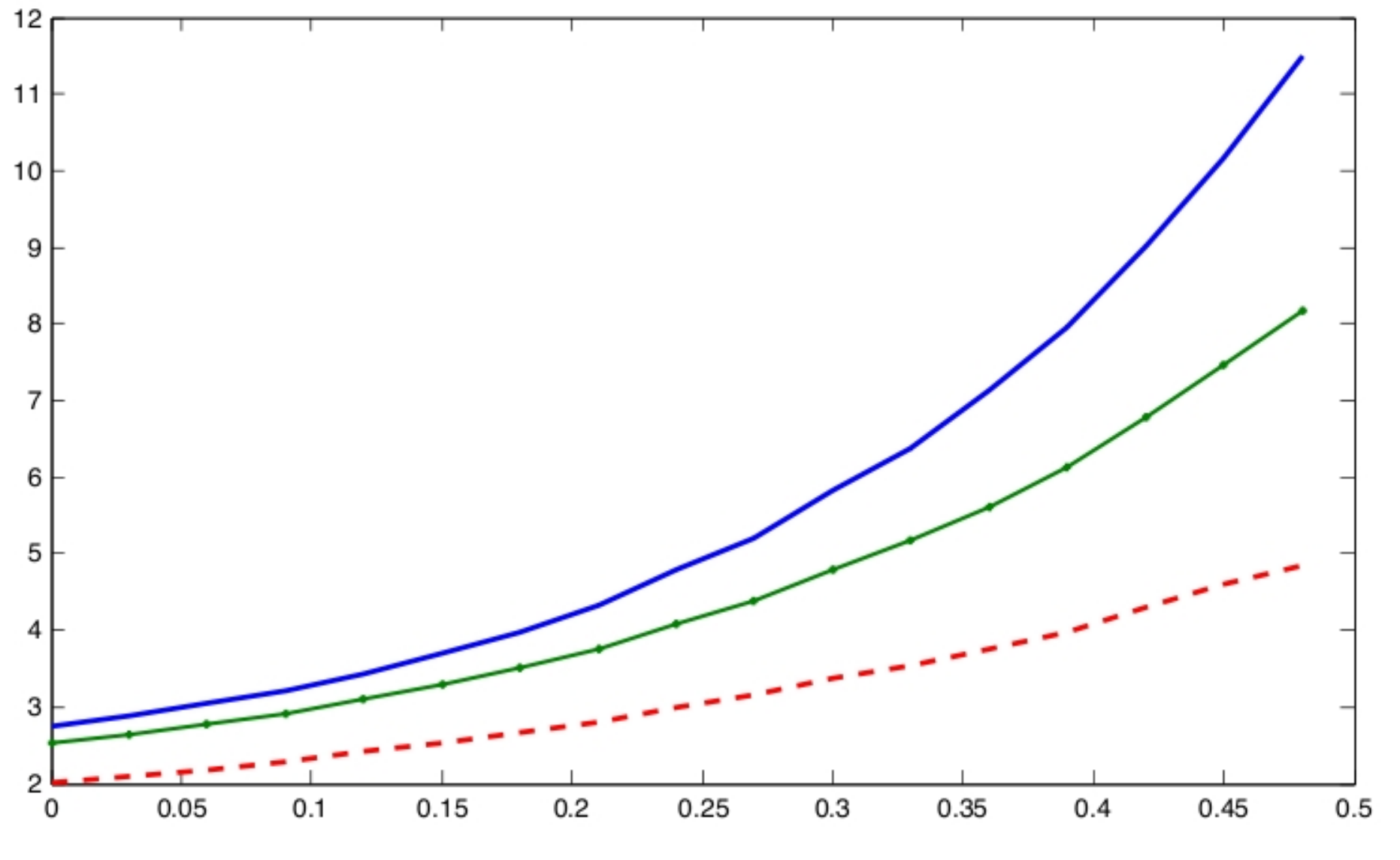}}
  \caption{The upper $5\%$ quantile of $I(B^{\rm I}_{d+.5})$ as a function of $d$,
  for $\underline\tau=.05$ (solid graph);
  $\underline\tau=.1$ (solid graph with points); $\underline\tau=.2$ (dashed graph).
  The right plot is a zoom in of the left plot in the region $d\in [0,.5)$.}
  \label{fig:qq}
\end{figure}

\subsection{Empirical comparisons of the tests $I_n$, $R_n$, $W_n$, $I^{Kim}_n$, $R^{Kim}_n$, and $W^{Kim}_n$}
\label{sec:comp-test-stat}

In this section we compare the test procedures based on our
statistics $I_n$, $R_n$, $W_n$  and the corresponding Kim's statistics
\begin{equation}\label{kim_integral} I_n^{Kim} := \int_{\mathcal{T}} {\cal K}_n(\tau) \d \tau, \quad
W_n^{Kim} := \sup_{\tau\in\mathcal{T}} {\cal K}_n(\tau) \quad \text{and}\quad   R_n^{Kim} :=  \frac{\inf_{\tau\in {\cal T}}
U^*_{n-\lfloor n\tau\rfloor}(X)}{\inf_{\tau\in {\cal T}}U_{\lfloor n\tau\rfloor}(X)},
\end{equation}
where ${\cal K}_n(\tau)$,  $U_k(X)$ and $U^*_{n-k}(X)$ are  defined in \eqref{kim0} and \eqref{kim00}.

The empirical size of the above tests is evaluated by simulating the
FARIMA(0,$d$,0) process of Example \ref{ex1}  for $d = 0, \, 0.1, \, 0.2, \, 0.3, \, 0.4$. The empirical power is estimated
by simulating the FARIMA process of Example \ref{ex1}, (\ref{changed}) with
of $d_1, d_2 \in \{0, \, 0.1, \, 0.2, \, 0.3, \, 0.4\}$ and the change-point of $d$ in the middle of the sample $(\theta^* = 0.5)$.

Tables \ref{compa500} and \ref{compa5000}
display the estimated level and power based on $10^4$ replications of
the testing procedures, for respective sample sizes  $n=500$ and $n=5000$.
These results show that for all six statistics and three values of $\underline\tau = 0.05, 0.1, 0.2$, the
estimated level is close to the nominal level $5\%$.
We also observe that while the performance of the tests $I_n$ and $I_n^{Kim}$ does not much depend
on $\underline\tau$, the last property is not shared by $W_n^{Kim}$, $R_n^{Kim}$, $W_n$ and $R_n$ .

Tables \ref{compa500} and \ref{compa5000} suggest that when $\underline\tau$ is small, $I_n$ clearly outperforms the remaining five tests.
As $\underline\tau$ increases, the performance of $R_n$ becomes comparable
to that of $I_n$, while $W_n$, $I_n^{Kim}$, $W_n^{Kim}$ and $R_n^{Kim}$ still remain less efficient.
Clearly, it make sense to choose $\underline\tau $ as small as possible, since none of these tests can detect
a change-point that occurs outside the testing interval $\mathcal{T} =[\underline\tau , 1-
  \underline\tau]$.


In conclusion, given the choice $\underline\tau = 0.05$ (or $\underline\tau = 0.1$), the statistic
  $I_n$ seems preferable to $W_n$ and $R_n$ and the three
Kim's statistics $I^{Kim}_n$, $R^{Kim}_n$, and $W^{Kim}_n$.

\subsection{Further simulation results pertaining to the test $I_n$}

As noted above, the results in Subsection 7.1 suggest that $I_n$ is favorable among the six tests
in the case when the observations follow a `pure' FARIMA(0,$d$,0) model
 with a possible jump of $d$ for $d \in [0, .5)$.
Here, we explore the performance of $I_n$ when the observations are simulated from other
classes
of processes following the hypotheses $\bf{H_0[I]}$ and $\bf{H_1[I]}$.

Table~\ref{fd500} extends the results of Tables~\ref{compa500} and \ref{compa5000}
to a larger interval of $d$ values, viz., $0<d_1\leq d_2<1.5$. Recall that, in 
accordance with Definition~\ref{defId} (iii), for $d>.5$ the FARIMA$(0,d,0)$ process is defined by
\begin{equation}\label{Fdarb}
    X_t = \sum_{i=1}^t Y_i, \qquad t=1, \dots, n,
\end{equation}
 where $\{Y_t\}$ is a stationary FARIMA$(0,d-1,0)$.

Figure~\ref{fig:traj} shows some trajectories simulated from
model (\ref{changed}) with fixed $d_2-d_1 =0.3$ and three different values of $d_1$, with the change-point in the middle
of the sample. From visual inspection of these paths, it seems that it is more difficult to
detect a change in the memory parameter when 
$0\leq d_1 < d_2 < .5$  or
$.5< d_1 < d_2$ (top or bottom graphs) than when $d_1< .5 <  d_2$ (middle graph). Note that the top and bottom graphs of  Figure~\ref{fig:traj}
correspond to $d_1, d_2$ belonging to the same `stationarity interval'
(either $[0, .5)$ or $(.5, 1.5)$) and the middle graph to $d_1, d_2$  falling into different `stationarity intervals'  $[0, .5)$ and $(.5, 1.5)$.
The above visual observation is indeed confirmed in  Table~\ref{fd500}.  The last table also shows that when
the difference $d_2-d_1$ is fixed, the test $I_n$
is more powerful in the region  $0\leq  d_1 < d_2<.5$ than in
$.5 < d_1 < d_2$.

Tables~\ref{fd500arp} and \ref{fd500arm} illustrate the performance of $I_n$ when a
positive, resp.\ negative, autoregressive part is added to the fractional process $\{\varepsilon_t(d_i)\}$ ($i=1,2$) in the
model \eqref{changed}.
These tables 
show that the  performance of the test is essentially preserved, especially when the autoregressive coefficient is positive.
However in the case of negative  autoregressive coefficient
the estimated level is slightly more disturbed (Table~\ref{fd500arm}).
Tables~\ref{fd500}, \ref{fd500arp} and \ref{fd500arm} also confirm that $I_n$ is not
very sensitive to the choice of the parameter $\underline\tau$, i.e.\ to the length of testing interval.

Finally, we assess the power of the test  for fractionally integrated models with changing memory parameters discussed in Section~\ref{sec:fract-integr-models}. Figure~\ref{fig:trajtv} presents sample paths of the {\it rapidly changing memory} model
in \eqref{X1t} and the {\it gradually changing memory} model in \eqref{X2t},
for the same function $d(t/n)  = .2 + .6\; t/n$ with the middle point $t = \lfloor n/2\rfloor$ marking the transition from
`stationarity regime' $d \in [0, .5)$ to `nonstationarity regime' $d \in (.5, 1)$. The visual impression
from Figure~\ref{fig:trajtv} is that the above transition is much easier to detect
for the rapidly changing memory model than for the gradually changing memory model.

Table~\ref{fd500trend} displays the estimated power of the test $I_n$  for the rapidly changing memory model \eqref{X1t}
when $d(\tau) = d_1+(d_2-d_1)\tau$, $\underline\theta =0$ and $\overline\theta =1$.
The null hypothesis is  naturally less often rejected for this model than for the model defined in (\ref{changed}),
cf.\ Table~\ref{fd500}. However, the estimated power still seems to be satisfactory.
Similar simulations under gradually changing memory model (not included in this paper) show that the test has more difficulty to
detect this type of changing memory on small samples. However, when  $n$ is larger than 500, the difference in the estimated
power between the gradually and rapidly memory cases becomes negligible.

\subsection{Simulations  in the presence of linear trend}

In this section we illustrate the performances of the test based on the
de-trended statistic $\mathcal{I}_n$ defined in \eqref{trstat1}.
This  testing  procedure is implemented similarly  to the previous
one.  Note that  the critical region still  depends on  the
memory  parameter $d$ which is estimated as follows:
having observed $X_1,\dots,X_n$ we estimate $d$ using NELWE estimate on
  the residuals from the least-squares
regression of $(X_j)_{1\le j \le n} $ on $(a + b j)_{1\le j \le n}$.

First we apply the de-trended test on series without trend, namely from model
\eqref{changed}. Table \ref{trend} displays the estimated level and power of
the test for this model.  The  estimated level is close to the nominal level.
Moreover the power is close to that obtained in Table~\ref{fd500}. Therefore
the performances of the testing procedure are  preserved
even if the estimation of the linear trend was not necessary.

Second we assess the power of this test in presence of a linear
trend (see Table~\ref{trend2}). Figure~\ref{fig:trend} presents sample paths of  models
defined in \eqref{eq:Xtrend} with $a=1$, $b=.01$, $\theta^*=1/2$,
$n=500$ and different values of $d_1$, $d_2$. For this model, Table~\ref{trend2}
summarizes the estimated level and power, which are similar to Table~\ref{trend}.

\begin{figure}[h]
  \centering
  \includegraphics[width=\textwidth , height=.55\textheight]{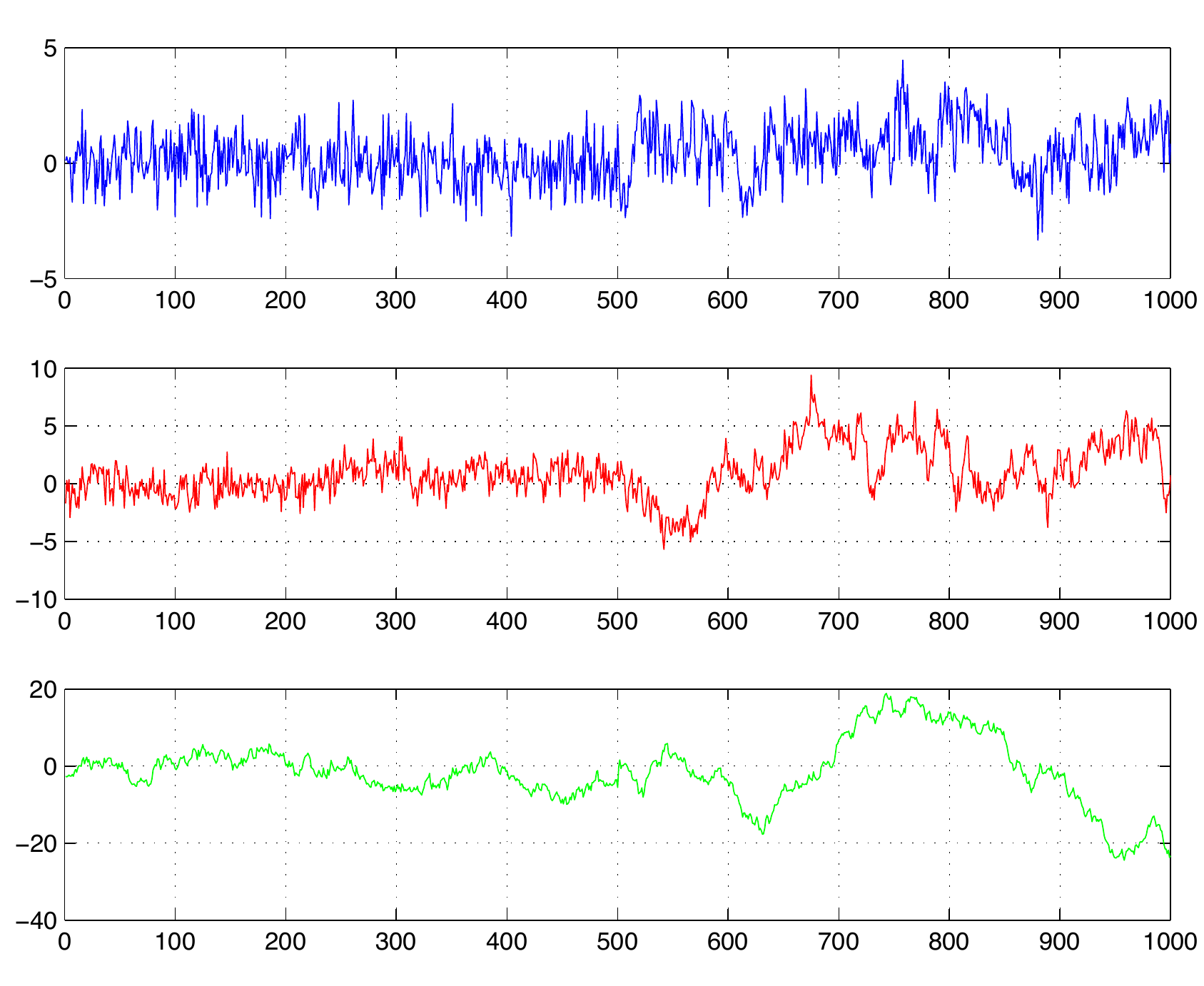}
  \caption{Sample paths simulated from model (\ref{changed}) with $\theta^*=1/2$ for different values of $d_1$ and $d_2$:
$d_1=.1$, $d_2=.4$ (top); $d_1=.3$, $d_2=.6$ (middle); $d_1=.8$, $d_2=1.1$ (bottom). The sample size is $n=1000$. }
\label{fig:traj} \end{figure}

\begin{figure}[h]
  \centering
 \includegraphics[width=\textwidth, height=.5\textheight]{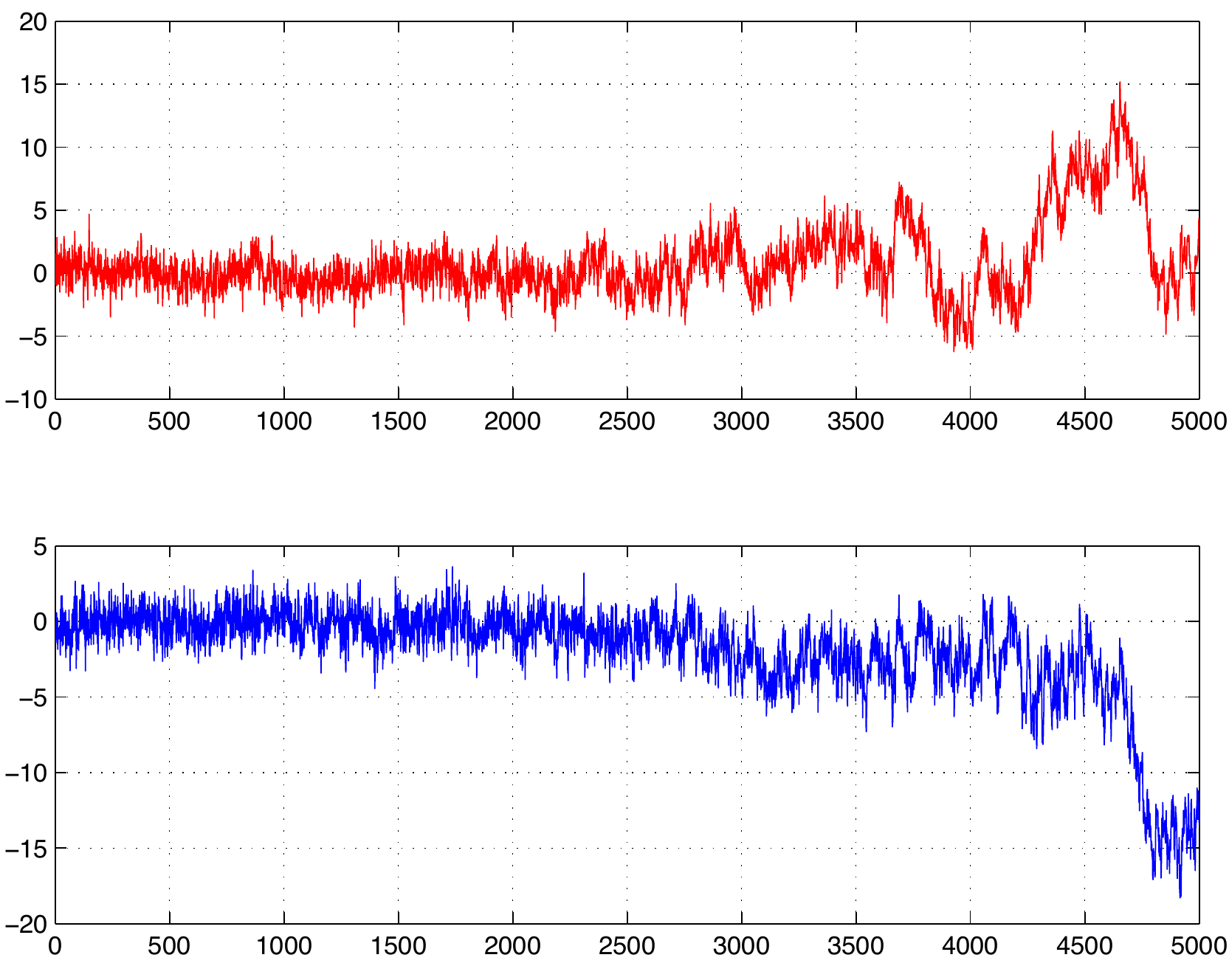}
 \caption{Sample paths simulated from models \eqref{X1t} (top)
 and \eqref{X2t} (bottom) with $d(\tau) = .2 + .6 \tau$, $\underline\theta=0$ and $\overline\theta=1$.}
\label{fig:trajtv}
\end{figure}

\begin{figure}[h]
  \centering
  \includegraphics[width=\textwidth, height=.4\textheight]{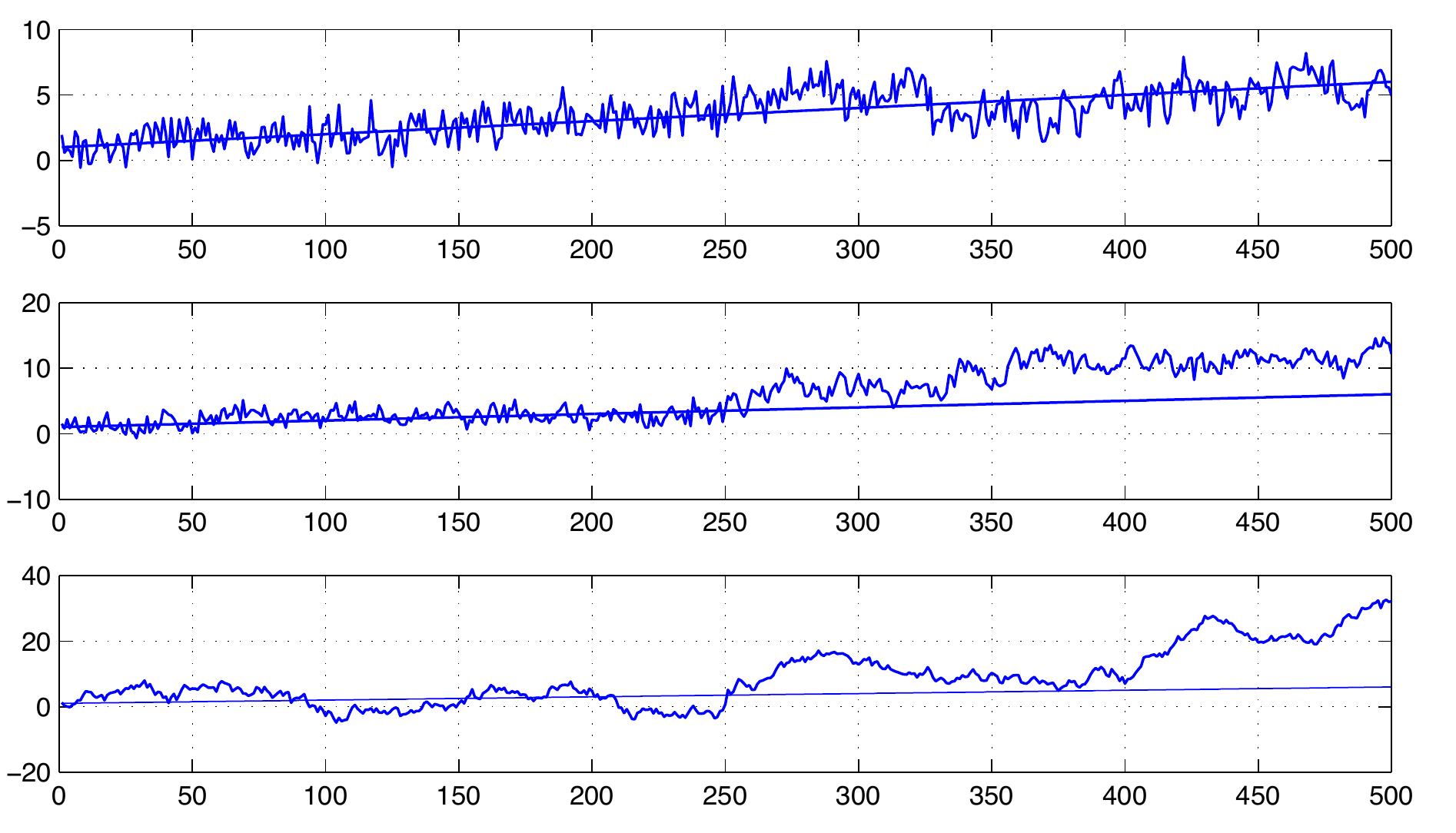}
\caption{Comparison of trajectories simulated from the model \eqref{eq:Xtrend}
with $a=1$, $b=.01$, change point $\theta^*=1/2$ and different values of $d_1$ and $d_2$:
$d_1=.1$, $d_2=.4$ (top); $d_1=.3$, $d_2=.6$ (middle); $d_1=.8$, $d_2=1.1$ (bottom).
The sample size is $n=500$.}
\label{fig:trend} \end{figure}

\clearpage

\begin{table}
  \centering
 \begin{tabular}{l|lllll|lllll|lllll}
\hline\hline
\multicolumn{1}{c}{} &
\multicolumn{5}{c}{$\underline\tau=0.05$}&\multicolumn{5}{c}{$\underline\tau=0.1$}&\multicolumn{5}{c}{$\underline\tau=0.2$}\\
\hline
 \backslashbox{$d_2$}{$d_1$}&0&0.1&0.2&0.3&0.4&0&0.1&0.2&0.3&0.4&0&0.1&0.2&0.3&0.4\\ \hline\hline
\multicolumn{16}{c}{$W_n$ statistic}\\ \hline
0& 3.6&&&&&2.4&&&&&2.3&&&&\\
0.1&10.0&4.4&&&&6.3&2.5&&&&10.4&2.6&&&\\
0.2&19.0&10.2&4.4&&&11.8&5.1&2.8&&&25.2&9.9&2.7&&\\
0.3&25.9&16.9&9.5&4.9&&15.9&9.0&4.7&3.1&&45.1&23.3&9.4&2.9&\\
0.4&31.0&21.6&13.1&8.3&4.8&19.9&12.2&6.8&4.5&2.8&65.1&42&21.0&8.7&2.6\\\hline
\multicolumn{16}{c}{$R_n$ statistic}\\ \hline
0& 3.4&&&&&2.6&&&&&2.5&&&&\\
0.1&11.2&3.9&&&&8.0&2.9&&&&11.6&2.9&&&\\
0.2&24.5&12.0&4.6&&&17.0&7.2&3.1&&&28.5&11.3&3.2&&\\
0.3&37.6&23.1&11.5&4.9&&27.2&14.7&6.6&3.4&&50.6&27.3&11.1&3.4&\\
0.4&49.5&34.0&19.8&10.6&5.0&41.1&23.7&12.4&6.7&3.0&72.0&49.0&25.5&10.5&3.4\\\hline
\multicolumn{16}{c}{$I_n$ statistic}\\ \hline
0&2.3&&&&&2.5&&&&&2.9&&&&\\
0.1&10.4&2.6&&&&11.6&2.9&&&&13.2&3.6&&&\\
0.2&25.2&9.9&2.7&&&28.5&11.3&3.2&&&31.8&13.4&3.8&&\\
0.3&45.1&23.3&9.4&2.9&&50.6&27.3
&11.1&3.4&&55.0&31.0&13.0&4.0&\\
0.4&65.1&42.0&21.0&8.7&2.6&72.0&49.0&25.5&10.5&3.4&77.1&54.3&29.6&12.4&4.1\\
\hline
\multicolumn{16}{c}{$W_n^{Kim}$ statistic}\\ \hline
0& 3.1&&&&&3.3&&&&&3.3&&&&\\
0.1&6.6&3.2&&&&8.3&3.7&&&&10.3&4.2&&&\\
0.2&10.7&6.2&3.3&&&14.6&7.4&3.6&&&20.2&9.5&4.4&&\\
0.3&16.5&9.0&5.3&3.1&&25.2&13.6&7.0&3.8&&36.6&20&9.4&4.4&\\
0.4&24.4&15.5&9.5&5.7&3.4&41.7&26.4&15.7&7.7&4.2&58.3&37.9&21.7&10.0&4.8\\\hline
\multicolumn{16}{c}{$R_n^{Kim}$ statistic}\\ \hline
0& 3.7&&&&&3.3&&&&&3.3&&&&\\
0.1&9.4&4.3&&&&9.9&4.4&&&&10.7&4.3&&&\\
0.2&17.1&9.5&4.7&&&21.1&11.1&4.7&&&24.1&11.6&4.5&&\\
0.3&25.8&15.3&8.8&4.6&&34.9&20.6&10.6&5.1&&42.4&24.2&11.7&4.8&\\
0.4&31.4&22.1&13.6&8.4&4.9&47.4&32.3&19.1&10.3&5.0&58.8&40.2&22.5&11.7&4.9\\\hline
\multicolumn{16}{c}{$I_n^{Kim}$ statistic}\\ \hline
0& 2.6&&&&&2.8&&&&&3.1&&&&\\
0.1&8.6&3.2&&&&9.5&3.5&&&&10.5&4.0&&&\\
0.2&18.7&8.3&3.3&&&20.8&9.1&3.5&&&23.1&10.4&4.0&&\\
0.3&35.0&18.3&8.4&3.1&&38.9&21.0&9.3&3.4&&42.0&23.4&10.6&4.0&\\
0.4&57.0&37.3&20.5&9.0&3.9&62.3&41.3&23.5&10.2&4.2&64.8&43.8&24.5&11.3&4.4\\\hline
\hline
\end{tabular}
\caption{Estimated level ($d_1=d_2$) and power ($d_1\not=d_2$) (in \%) of the tests
$I_n$, $W_n$, $R_n$, $I^{Kim}_n$, $W^{Kim}_n$, $R^{Kim}_n$.
The nominal level is $\alpha = 5\%$. The samples are simulated from model (\ref{changed}) with $\theta^*=1/2$.
The sample size is $500$ and the number of independent replications is $10^4$.}\label{compa500}
\end{table}

\clearpage

\begin{table}
  \centering
 \begin{tabular}{l|lllll|lllll|lllll}
\hline\hline
\multicolumn{1}{c}{} &
\multicolumn{5}{c}{$\underline\tau=0.05$}&\multicolumn{5}{c}{$\underline\tau=0.1$}&\multicolumn{5}{c}{$\underline\tau=0.2$}\\
\hline
 \backslashbox{$d_2$}{$d_1$}&0&0.1&0.2&0.3&0.4&0&0.1&0.2&0.3&0.4&0&0.1&0.2&0.3&0.4\\ \hline\hline
\multicolumn{16}{c}{$W_n$ statistic}\\ \hline
0&3.4&&&&&3.6&&&&&4.1&&&&\\
0.1&20.2&3.7&&&&24.2&4.1&&&&28.3&4.5&&&\\
0.2&50.1&16.0&3.8&&&61.1&20.9&4.5&&&68.6&26.0&4.7&&\\
0.3&74.6&38.5&13.3&3.6&&87.1&51.1&17.7&4.1&&93.1&62.0&23.6&4.7&\\
0.4&88.7&64.1&31.2&11.6&3.7&96.6&79.6&44.3&15.0&4.0&99.2&88.8&56.7&20.0&4.4\\ \hline
\multicolumn{16}{c}{$R_n$ statistic}\\ \hline
0&3.7&&&&&3.9&&&&&4.1&&&&\\
0.1&29.0&5.1&&&&30.9&4.9&&&&31.3&5.0&&&\\
0.2&65.0&25.3&4.9&&&71.0&28.1&4.8&&&73.1&30.5&4.6&&\\
0.3&85.3&53.8&21.2&5.0&&91.7&62.2&24.4&4.9&&94.3&67.0&26.5&4.9&\\
0.4&94.1&75.7&44.8&18.4&4.8&98.0&84.8&54.1&21.7&4.6&99.3&90.3&60.3&24.4&4.5\\ \hline
\multicolumn{16}{c}{$I_n$ statistic}\\ \hline
0&2.9&&&&&3.2&&&&&3.6&&&&\\
0.1&28.5&3.4&&&&29.5&3.8&&&&30.7&4.1&&&\\
0.2&73.4&27.8&3.6&&&74.0&28.8&3.9&&&73.7&30.1&4.2&&\\
0.3&95.9&68.0&24.8&3.5&&96.3&69.8&26.5&3.7&&95.9&69.7&27.9&4.3&\\
0.4&99.5&92.3&62.5&21.8&3.5&99.7&93.6&65.7&24.1&3.8&99.7&94.1&66.6&25.6&4.0\\ \hline
\multicolumn{16}{c}{$W_n^{Kim}$ statistic}\\ \hline
0&5.0&&&&&4.9&&&&&4.8&&&&\\
0.1&20.3&5.1&&&&22.0&5.1&&&&23.8&5.2&&&\\
0.2&44.9&15.7&4.7&&&51.1&18.1&5.1&&&55.2&21.2&5.5&&\\
0.3&68.8&37.2&13.4&4.6&&77.6&45.7&16.6&5.0&&83.1&52.1&19.5&5.3&\\
0.4&83.7&60.1&32.4&12.8&4.2&92.2&72.4&41.4&14.9&4.5&96.0&79.7&48.1&18.2&4.6\\
\hline
\multicolumn{16}{c}{$R_n^{Kim}$ statistic}\\ \hline
0&3.9&&&&&4.2&&&&&4.5&&&&\\
0.1&26.5&4.8&&&&26.5&4.6&&&&26.6&5.1&&&\\
0.2&59.6&23.1&4.9&&&63.2&24.9&5.1&&&62.7&24.6&5.1&&\\
0.3&81.6&49.9&19.6&5.4&&87.0&56.4&22.4&5.2&&87.4&57.4&22.6&5.0&\\
0.4&91.9&72.1&42.0&17.6&4.9&96.5&80.8&50.1&20.4&4.9&97.2&82.6&52.5&20.1&4.5\\
\hline
\multicolumn{16}{c}{$I_n^{Kim}$ statistic}\\ \hline
0&3.5&&&&&3.7&&&&&4.2&&&&\\
0.1&23.5&3.6&&&&24.1&4.0&&&&24.6&4.3&&&\\
0.2&58.4&21.4&4.3&&&58.5&22.0&4.4&&&57.4&22.4&4.7&&\\
0.3&86.1&54.8&19.9&4.0&&86.8&55.8&20.5&4.2&&86.1&55.4&20.8&4.5&\\
0.4&96.6&82.2&52.1&18.4&4.0&97.3&83.6&53.6&19.2&3.7&97.3&83.4&53.4&19.7&4.0\\
\hline\hline
\end{tabular}
\caption{Estimated level ($d_1=d_2$) and power ($d_1\not=d_2$) (in \%) of the tests
$I_n$, $W_n$, $R_n$,  $I^{Kim}_n$, $W^{Kim}_n$, $R^{Kim}_n$.
The nominal level is $\alpha = 5\%$. The samples are simulated from model (\ref{changed}) with $\theta^*=1/2$.
The sample size is $5000$ and the number of independent replications is $10^4$.}\label{compa5000}
\end{table}

\clearpage

\begin{sidewaystable}
      \centering
      \footnotesize
\begin{tabular}{l|l|lllll;{.4pt/1pt}lllll|lllll;{.4pt/1pt}lllll}
\hline\hline
&\multicolumn{1}{c}{} &\multicolumn{10}{c|}{$n=500$}& \multicolumn{10}{c}{$n=5000$}\\
\cline{2-22}
$\underline\tau$ &\backslashbox{$d_2$}{$d_1$}&0&0.1&0.2&0.3&0.4&0.6&0.8&1&1.2&1.4&0&0.1&0.2&0.3&0.4&0.6&0.8&1&1.2&1.4\\
\hline\hline
0.05&0&2.3   &&&&&&&&&&2.9&&&&&&&&&\\
&0.1&10.4&2.6    &&&&&&&&&28.5&3.4&&&&&&&&\\
&0.2&25.2&9.9&2.7     &&&&&&&&73.4&27.8&3.6&&&&&&&\\
&0.3&45.1&23.3&9.4&2.9 &&&&&&&95.9&68.0&24.8&3.5&&&&&&\\
&0.4&65.1&42.0&21.0&8.7&2.6          &&&&&&99.5&92.3&62.5&21.8&3.5&&&&&\\ \cdashline{2-22}[.4pt/1pt]
&0.6&92.7&85.4&71.8&54.5&35.2&4.9    &&&&&100&99.9&99.7&97.1&88.0&5.3&&&&\\
&0.8&94.7&88.8&78.8&63.1&45.2&7.8&6.2    &&&&100&100&99.8&98.7&92.6&11.3&5.0&&&\\
&1&97.2&94.4&89.2&80.8&67.4&19.8&12.5&6.7&&&100&100&100&99.9&98.9&45.3&28.1&5.1&&\\
&1.2&98.6&97.6&95.5&91.7&86.9&55.5&44.9&23.5&7.2&&100&100&100&100&99.9&87.4&78.6&41.4&5.3&\\
&1.4&99.5&99.1&98.1&97.4&95.4&84.3&78.8&67.6&44.2&4.9&100&100&100&100&100&98.0&96.4&85.9&59.0&3.3\\\hline
0.1&0&2.5 &&&&&&&&&&3.2&&&&&&&&&\\
&0.1&11.6&2.9&&&&&&&&&29.5&3.8&&&&&&&&\\
&0.2&28.5&11.3&3.2 &&&&&&&&74.0&28.8&3.9&&&&&&&\\
&0.3&50.6&27.3&11.1&3.4&&&&&&&96.3&69.8&26.5&3.7&&&&&&\\
&0.4&72.0&49.0&25.5&10.5&3.4 &&&&&&99.7&93.6&65.7&24.1&3.8&&&&&\\ \cdashline{2-22}[.4pt/1pt]
&0.6&97.9&94.0&84.5&68.1&47.3&4.9 &&&&&100&100&99.9&98.7&92.5&5.2&&&&\\
&0.8&98.8&96.3&90.9&78.4&59.8&9.4&5.9 &&&&100&100&99.9&99.4&96.0&11.8&4.9&&&\\
&1&99.4&98.8&96.7&92.3&82.4&27.5&16.8&6.2&&&100&100&100&100&99.7&51.9&34.0&5.1&&\\
&1.2&99.8&99.6&99.0&97.6&95.0&65.1&52.5&26.0&6.3&&100&100&100&100&100&91.7&84.2&44.6&5.1&\\
&1.4&100&99.8&99.7&99.3&98.8&89.7&84.2&70.3&42.4&3.7&100&100&100&100&100&99.1&98.2&88.9&57.8&2.9\\ \hline
0.2&0&2.9 &&&&&&&&&&3.6&&&&&&&&&\\
&0.1&13.2&3.6&&&&&&&&&30.7&4.1&&&&&&&&\\
&0.2&31.8&13.4&3.8&&&&&&&&73.7&30.1&4.2&&&&&&&\\
&0.3&55.0&31.0&13&4.0&&&&&&&95.9&69.7&27.9&4.3&&&&&&\\
&0.4&77.1&54.3&29.6&12.4&4.1&&&&&&99.7&94.1&66.6&25.6&4.0&&&&&\\ \cdashline{2-22}[.4pt/1pt]
&0.6&99.3&97.6&91.4&78.3&57.2&5.2&&&&&100&100&99.9&99.4&95.4&5.1&&&&\\
&0.8&99.7&98.9&95.8&87.7&71.8&10.6&5.4 &&&&100&100&100&99.8&98.1&11.6&4.9&&&\\
&1&99.9&99.8&99.2&97.3&91.7&34.7&21.0&5.8&&&100&100&100&100&99.9&57.4&39.4&4.6&&\\
&1.2&100&99.9&99.8&99.6&98.7&74.1&60.5&28.6&5.2&&100&100&100&100&100&94.8&89.5&48.4&4.9&\\
&1.4&100&100&99.9&99.9&99.8&94.5&89.6&73.4&41.4&2.7 &100&100&100&100&100&99.7&99.3&91.8&57.8&2.7\\ \hline\hline
\end{tabular}
\caption{Estimated level ($d_1=d_2$) and power ($d_1\not=d_2$) (in \%) of the test  $I_n$
 with nominal level $\alpha=5\%$. The samples are simulated from model (\ref{changed})
 with $\theta^*=1/2$. The sample sizes are $n=500$ and $n=5000$, the number of independent replications is  $10^4$. }\label{fd500}
\end{sidewaystable}
\clearpage

\begin{sidewaystable}
      \centering \footnotesize
\begin{tabular}{l|l|lllll;{.4pt/1pt}lllll|lllll;{.4pt/1pt}lllll}
\hline\hline
&\multicolumn{1}{c}{} &\multicolumn{10}{c|}{$n=500$}& \multicolumn{10}{c}{$n=5000$}\\
\cline{2-22}
$\underline\tau$ &\backslashbox{$d_2$}{$d_1$}&0&0.1&0.2&0.3&0.4&0.6&0.8&1&1.2&1.4&0&0.1&0.2&0.3&0.4&0.6&0.8&1&1.2&1.4\\
\hline\hline
0.05&0&2.6&&&&&&&&& &3.0&&&&&&&&&\\
&0.1&10.7&3&&&&&&&& &28.6&3.6&&&&&&&&\\
&0.2&27.2&10.8&3.4&&&&&&& &73.4&26.9&3.6&&&&&&&\\
&0.3&47.3&24.0&9.1&2.9&&&&&& &95.2&68.3&23.7&3.7&&&&&&\\
&0.4&66.4&42.6&22.1&8.1&2.9&&&&& &99.4&91.5&63.5&21.5&3.4&&&&&\\ \cdashline{2-22}[.4pt/1pt]
&0.6&92.6&84.4&70.6&53.7&35.0&5.0&&&& &100&100&99.6&97.4&87.5&5.2&&&&\\
&0.8&94.6&89.2&79.1&63.5&43.9&8.3&5.9&&& &100&99.9&99.9&98.7&92.2&11.3&5.4&&&\\
&1&97.3&94.5&89.0&80.7&67.3&20.2&13.1&6.9&& &100&100&99.9&99.9&99.0&44.8&28.1&5.5&&\\
&1.2&98.5&97.7&95.4&92.0&86.8&55.4&44.3&23.2&7.1& &100&100&100&99.9&99.9&87.3&78.0&43.0&5.7&\\
&1.4&99.5&99.1&98.4&97.0&95.1&83.4&78.9&67.6&44.7&4.8 &100&100&100&100&100&98.2&96.3&84.7&58.9&3.1\\ \hline
0.1&0&2.9&&&&&&&&& &3.3&&&&&&&&&\\
&0.1&11.9&3.3&&&&&&&& &29.5&4.0&&&&&&&&\\
&0.2&29.7&12.2&3.8&&&&&&& &73.6&27.8&3.8&&&&&&&\\
&0.3&52.2&27.4&10.7&3.3&&&&&& &95.7&70.0&24.9&4.0&&&&&&\\
&0.4&72.7&48.5&26.1&9.9&3.3&&&&& &99.6&92.9&66.2&23.5&4.0&&&&&\\ \cdashline{2-22}[.4pt/1pt]
&0.6&97.7&92.9&83.6&67.7&47.2&5.0&&&& &100&100&99.9&98.8&92.1&5.1&&&&\\
&0.8&98.9&96.5&90.6&78.0&58.0&9.4&5.6&&& &100&99.9&99.9&99.5&96.1&11.9&5.2&&&\\
&1&99.4&98.7&96.9&92.0&82.5&26.9&16.6&6.5&& &100&100&100&100&99.7&51.0&34.9&5.4&&\\
&1.2&99.8&99.6&98.9&97.5&94.8&65.4&52.1&25.7&6.2& &100&100&100&100&99.9&91.6&84.0&46.0&5.4&\\
&1.4&99.9&99.9&99.7&99.3&98.6&89.3&84.2&70.1&43.2&3.9 &100&100&100&100&100&99.1&98.2&87.9&57.8&2.8\\ \hline
0.2&0&3.2&&&&&&&&& &3.7&&&&&&&&&\\
&0.1&13.5&4.0&&&&&&&& &30.6&4.5&&&&&&&&\\
&0.2&32.8&13.8&4.3&&&&&&& &73.4&29.0&4.3&&&&&&&\\
&0.3&56.5&30.8&12.6&4.0&&&&&& &95.4&70.1&26.6&4.3&&&&&&\\
&0.4&77.4&53.6&29.5&12.2&4.2&&&&& &99.7&93.2&67.4&25.1&4.3&&&&&\\ \cdashline{2-22}[.4pt/1pt]
&0.6&99.4&96.9&90.6&77.6&57.0&5.1&&&& &100&100&99.9&99.4&94.6&5.0&&&&\\
&0.8&99.8&98.9&96&87.4&70.0&10.4&5.3&&& &100&100&100&99.8&98.0&11.7&5.1&&&\\
&1&99.9&99.8&99.3&97.0&91.8&33.6&20.3&5.8&& &100&100&100&100&99.9&56.6&40.5&5.0&&\\
&1.2&99.9&99.9&99.9&99.6&98.6&74.8&60.1&28.6&5.8& &100&100&100&100&100&94.7&89.4&49.3&4.9&\\
&1.4&100&100&100&99.9&99.7&94.1&89.3&73.2&41.4&3.1 &100&100&100&100&100&99.7&99.4&91.2&56.9&2.2\\ \hline\hline
\end{tabular}
\caption{Estimated level ($d_1=d_2$) and power ($d_1\not=d_2$) (in \%) of the test $I_n $
  with nominal level $\alpha = 5\%$. The samples are simulated from \eqref{changed} with the FARIMA(0,$d$,0) processes
  replaced by FARIMA$(1,d,0)$ models with the AR parameter $.7$ and $\theta^*=1/2$.
The sample sizes are $n=500$ and $n=5000$, the number of independent replications is $10^4$.} \label{fd500arp}
\end{sidewaystable}

\clearpage

\begin{sidewaystable}
      \centering \footnotesize
\begin{tabular}{l|l|lllll;{.4pt/1pt}lllll|lllll;{.4pt/1pt}lllll}
\hline\hline
&\multicolumn{1}{c}{} &\multicolumn{10}{c|}{$n=500$} & \multicolumn{10}{c}{$n=5000$}\\
\cline{2-22}
$\underline\tau$ &\backslashbox{$d_2$}{$d_1$}&0&0.1&0.2&0.3&0.4&0.6&0.8&1&1.2&1.4&0&0.1&0.2&0.3&0.4&0.6&0.8&1&1.2&1.4\\
\hline\hline
0.05&0&0.3&&&&&&&&& &1.8&&&&&&&&&\\
&0.1&2.4&0.6&&&&&&&& &24.3&2.7&&&&&&&&\\
&0.2&9.0&3.5&1.1&&&&&&& &70.2&25.3&3.1&&&&&&&\\
&0.3&23.3&11.6&4.7&1.3&&&&&& &94.6&66.5&22.7&3.3&&&&&&\\
&0.4&43.5&27.2&13.5&5.1&1.6&&&&& &99.5&92.7&63.1&22.5&3.7&&&&&\\\cdashline{2-22}[.4pt/1pt]
&0.6&92.9&87.6&77.4&64.2&46.5&3.5&&&& &99.9&100&99.7&97.9&89.3&6.5&&&&\\
&0.8&94.3&90.5&84.1&72.0&55.3&11.8&6.5&&& &100&99.9&99.8&98.9&93.0&12.6&5.3&&&\\
&1&96.7&94.9&91.2&84.4&74.3&25.9&16.8&8.0&& &100&100&99.9&99.8&99.1&46.3&29.2&5.8&&\\
&1.2&98.7&98.1&96.6&94.4&89.8&61.1&51.4&29.6&9.0& &100&100&100&100&99.9&87.8&78.8&42.8&6.0&\\
&1.4&99.5&99.2&98.6&98.1&96.6&87.2&83.0&71.6&51.1&7.5 &100&100&100&100&100&98.1&96.5&86.5&59.9&3.7\\\hline
0.1&0&0.1&&&&&&&&& &2.0&&&&&&&&&\\
&0.1&2.6&0.7&&&&&&&& &25.9&3.2&&&&&&&&\\
&0.2&11.0&4.5&1.4&&&&&&& &72.1&27.5&3.5&&&&&&&\\
&0.3&29.4&15.2&6.3&1.8&&&&&& &95.4&69.1&24.9&3.9&&&&&&\\
&0.4&55.0&37.2&19.2&7.8&2.5&&&&& &99.6&94.2&66.3&25.0&3.9&&&&&\\\cdashline{2-22}[.4pt/1pt]
&0.6&98.0&95.0&88.1&77.0&58.7&5.0&&&& &100&100&99.9&99.1&93.5&5.7&&&&\\
&0.8&98.9&97.0&93.8&85.1&69.8&14.9&7.5&&& &100&100&99.9&99.6&96.1&12.8&5.2&&&\\
&1&99.5&99.1&97.8&94.8&88.5&35.9&22.6&8.4&& &100&100&100&99.9&99.8&52.8&35.1&5.6&&\\
&1.2&99.9&99.6&99.4&98.6&96.8&72.7&60.7&32.7&8.8& &100&100&100&100&100&92.2&84.4&45.5&5.4&\\
&1.4&99.9&99.9&99.8&99.6&99.3&92.9&88.2&74.9&49.6&6.1 &100&100&100&100&100&99.3&98.2&89.2&58.3&3.1\\\hline
0.2&0&0.2&&&&&&&&& &2.5&&&&&&&&&\\
&0.1&3.7&1.0&&&&&&&& &27.8&3.9&&&&&&&&\\
&0.2&15.2&6.4&2.1&&&&&&& &72.9&29.2&4.1&&&&&&&\\
&0.3&37.8&20.8&9.4&2.8&&&&&& &95.4&70.4&26.6&4.5&&&&&&\\
&0.4&65.4&46.0&26.2&11.3&4.3&&&&& &99.7&94.5&67.6&26.4&4.3&&&&&\\\cdashline{2-22}[.4pt/1pt]
&0.6&99.4&98.0&93.4&84.0&66.3&5.9&&&& &100&100&100&99.6&95.7&5.4&&&&\\
&0.8&99.8&99.4&97.4&91.5&78.7&15.2&6.9&&& &100&100&99.9&99.9&97.9&12.1&5.5&&&\\
&1&99.9&99.9&99.6&98.6&95.0&42.2&26.0&7.3&& &100&100&100&100&99.9&57.2&41.3&5.3&&\\
&1.2&100&99.9&99.9&99.8&99.2&80.0&67.4&33.8&7.5& &100&100&100&100&100&95.3&89.9&49.3&4.8&\\
&1.4&100&100&99.9&99.9&99.9&96.2&92.4&77.3&47.2&4.3 &100&100&100&100&100&99.7&99.3&92.4&57.7&2.8\\
\hline\hline
\end{tabular}
\caption{Estimated level ($d_1=d_2$) and power ($d_1\not=d_2$) (in \%) of the test $I_n$
with nominal level $\alpha=5\%$. The samples are simulated from \eqref{changed} with the
FARIMA(0,$d$,0) processes replaced by FARIMA(1,$d$,0) models with the AR parameter $-.7$ and $\theta^*=1/2$.
The sample sizes are $n=500$ and $n=5000$, the number of independent replications is $10^4$.} \label{fd500arm}
\end{sidewaystable}

\clearpage

\begin{sidewaystable}
      \centering \footnotesize
\begin{tabular}{l|l|lllll;{.4pt/1pt}lllll|lllll;{.4pt/1pt}lllll}
\hline\hline
&\multicolumn{1}{c}{} &\multicolumn{10}{c|}{$n=500$} & \multicolumn{10}{c}{$n=5000$}\\
\cline{2-22}
$\underline\tau$ &\backslashbox{$d_2$}{$d_1$}&0&0.1&0.2&0.3&0.4&0.6&0.8&1&1.2&1.4&0&0.1&0.2&0.3&0.4&0.6&0.8&1&1.2&1.4\\
\hline\hline
0.05&0&2.4&&&&&&&&&   &2.8&&&&&&&&&\\
&0.1&6.4&2.7&&&&&&&&   &16.1&4.2&&&&&&&&\\
&0.2&13.9&6.8&3.1&&&&&&&   &42.3&15.2&3.4&&&&&&&\\
&0.3&23.3&12.2&6.3&3.0&&&&&&   &69.9&36.6&13.1&3.5&&&&&&\\
&0.4&33.1&20.4&10.7&5.2&2.9&&&&&   &86.8&61.2&31.8&11.2&3.2&&&&&\\ \cdashline{2-22}[.4pt/1pt]
&0.6&57.2&40.5&27.2&16.6&10.0&4.2&&&&   &97.8&91.1&74.6&49.7&26.5&4.9&&&&\\
&0.8&71.6&57.3&43.4&30.1&19.9&7.8&5.9&&&   &99.3&97.8&91.9&79.6&59.1&17.1&5.0&&&\\
&1&76.1&64.7&51.9&38.3&26.9&11.4&6.0&6.6&&   &99.8&99.1&97.4&92.1&82.8&41.6&11.1&5.3&&\\
&1.2&76.2&68.5&57&45.2&32.6&14.4&6.0&5.0&6.5&   &99.8&99.5&98.9&96.5&91.5&66.1&26.1&6.7&5.2&\\
&1.4&74.1&67.2&57.4&47.1&36.2&17.0&6.9&3.3&3.4&3.4   &99.8&99.6&99.2&97.7&94.5&78.3&40.4&11.2&3.3&2.3\\\hline
0.1&0&2.8&&&&&&&&&   &3.1&&&&&&&&&\\
&0.1&7.1&3.1&&&&&&&&   &16.3&4.2&&&&&&&&\\
&0.2&15.3&7.7&3.6&&&&&&&   &42.3&15.7&3.8&&&&&&&\\
&0.3&25.7&14.0&7.1&3.2&&&&&&   &69.6&36.5&13.8&3.9&&&&&&\\
&0.4&37.0&23.5&12.5&6.2&3.3&&&&&   &86.7&61.8&32.8&11.8&3.3&&&&&\\ \cdashline{2-22}[.4pt/1pt]
&0.6&63.8&46.6&31.8&19.9&11.9&4.5&&&&   &98.3&92.2&75.8&51.2&27.3&4.7&&&&\\
&0.8&80.5&66.9&52.0&36.5&24.3&9.5&5.7&&&   &99.8&98.8&94.0&82.7&62.3&18.5&4.9&&&\\
&1&87.7&78.7&65.4&50.8&36.4&15.2&7.4&6.2&&   &99.9&99.7&98.7&95.0&86.6&45.9&13.1&5.2&&\\
&1.2&90.8&85.2&74.8&62.2&47.7&22.1&8.8&6.3&5.8&   &99.9&99.9&99.7&98.2&94.7&72.1&29.9&8.4&4.9&\\
&1.4&90.9&86.0&78.5&67.7&55.1&27.7&10.9&4.9&4.1&2.8   &100&99.9&99.8&99.3&97.4&84.7&45.9&13.6&4.4&1.9\\ \hline
0.2&0&3.4&&&&&&&&&   &3.6&&&&&&&&&\\
&0.1&8.1&3.7&&&&&&&&   &16.3&4.8&&&&&&&&\\
&0.2&17.0&9.0&4.3&&&&&&&   &39.9&16.0&4.2&&&&&&&\\
&0.3&27.7&15.7&8.4&3.9&&&&&&   &65.8&35.4&14.0&4.3&&&&&&\\
&0.4&39.5&25.8&14.4&7.6&4.1&&&&&   &83.4&59.0&32.7&12.5&3.6&&&&&\\ \cdashline{2-22}[.4pt/1pt]
&0.6&65.6&49.8&34.5&22.1&13.6&4.6&&&&   &97.7&90.5&73.9&49.7&27.3&4.4&&&&\\
&0.8&83.1&70.9&56.3&40.4&26.9&10.6&4.9&&&   &99.7&98.6&93.1&81.5&60.5&19.9&5.0&&&\\
&1&92.0&83.5&71.3&57.2&42.2&18.4&8.6&5.9&&   &100&99.9&98.6&94.3&85.4&45.7&14.8&5.2&&\\
&1.2&95.9&91.3&82.5&70.5&55.6&27.0&11.6&7.0&4.9&   &99.9&99.9&99.7&98.4&94.3&70.6&30.2&9.9&4.8&\\
&1.4&97.1&93.6&87.3&77.1&64.6&33.3&14.4&6.4&5.1&2.7   &100&99.9&99.9&99.4&97.3&82.9&43.7&14.8&5.8&2.3\\
\hline\hline
\end{tabular}
\caption{Estimated level ($d_1=d_2$) and power ($d_1\ne d_2$) (in \%) of the test $I_n$ with nominal level
 $\alpha=5\%$. The samples are simulated from model \eqref{X1t} with $d(\tau)=d_1+(d_2-d_1)\tau$,
 $\underline\theta=0$ and $\overline\theta=1$. The sample sizes are $n=500$ and
$n=5000$, the number of independent replications is $10^4$.} \label{fd500trend}
\end{sidewaystable}

\begin{sidewaystable}
      \centering \footnotesize
\begin{tabular}{l|l|lllll;{.4pt/1pt}lllll|lllll;{.4pt/1pt}lllll}
\hline\hline
&\multicolumn{1}{c}{} &\multicolumn{10}{c}{$n=500$} & \multicolumn{10}{c}{$n=5000$}\\
\cline{2-22}
$\underline\tau$ &\backslashbox{$d_2$}{$d_1$}&0&0.1&0.2&0.3&0.4&0.6&0.8&1&1.2&1.4&0&0.1&0.2&0.3&0.4&0.6&0.8&1&1.2&1.4\\
\hline\hline
0.05&0&2.4&&&&&&&&&&3.0&&&&&&&&&\\
&0.1&11.5&3.3&&&&&&&&&32.8&3.6&&&&&&&&\\
&0.2&30.9&12.8&3.3&&&&&&&&79.2&31.5&4.0&&&&&&&\\
&0.3&53.4&30.5&12.9&4.3&&&&&&&97.2&72.9&28.5&4.3&&&&&&\\
&0.4&71.5&50.7&28.3&12.8&4.6&&&&&&99.7&94.0&67.6&25.5&4.4&&&&&\\\cdashline{2-22}[.4pt/1pt]
&0.6&88.7&80.7&68.2&52.7&36.4&6.2&&&&&100&99.9&99.7&97.2&87.2&5.6&&&&\\
&0.8&91.2&84.2&73.8&60.7&43.8&8.1&5.4&&&&100&100&99.9&98.8&92.6&11.1&4.6&&&\\
&1&93.9&88.9&81.8&71.3&57.8&16.6&10.2&4.0&&&100&100&100&99.9&99.4&47.0&27.0&4.7&&\\
&1.2&97.0&94.6&90.0&83.1&73.7&33.7&24.4&14.7&4.0&&100&100&100&100&99.9&86.2&73.8&32.7&4.3&\\
&1.4&99.5&98.9&97.6&95.2&90.8&62.0&50.7&34.1&31.9&4.9&100&100&100&100&100&98.5&95.5&75.6&47.6&4.3\\\hline
0.1&0&2.7&&&&&&&&&&3.3&&&&&&&&&\\
&0.1&12.6&3.6&&&&&&&&&33.7&3.9&&&&&&&&\\
&0.2&33.8&14.1&3.8&&&&&&&&79.8&32.6&4.2&&&&&&&\\
&0.3&57.8&33.8&14.1&4.7&&&&&&&97.3&74.6&30.1&4.4&&&&&&\\
&0.4&77.2&56.4&32.3&14.7&4.7&&&&&&99.8&95.4&71.2&27.4&4.6&&&&&\\
&0.6&96.4&91.2&81.9&65.7&46.0&5.5&&&&&100&100&99.9&98.9&92.2&5.4&&&&\\\cdashline{2-22}[.4pt/1pt]
&0.8&97.9&94.9&87.9&75.8&58.2&9.4&5.1&&&&100&100&100&99.7&96.6&12.5&4.6&&&\\
&1&99.3&98.0&95.0&89.2&78.3&23.5&13.2&4.3&&&100&100&100&99.9&99.9&55.4&33.1&5.0&&\\
&1.2&99.8&99.5&98.6&96.4&92.5&51.7&37.1&17.8&4.3&&100&100&100&100&100&92.8&84.2&39.0&4.4&\\
&1.4&99.9&99.9&99.8&99.6&98.6&82.1&71.0&45.9&33.2&4.7&100&100&100&100&100&99.7&98.7&85.5&49.5&4.5\\\hline
0.2&0&3.2&&&&&&&&&&3.8&&&&&&&&&\\
&0.1&14.1&4.4&&&&&&&&&34.0&4.3&&&&&&&&\\
&0.2&36.0&15.7&4.5&&&&&&&&79.2&32.7&4.6&&&&&&&\\
&0.3&61.0&37.0&15.6&5.3&&&&&&&97.0&74.7&30.7&4.9&&&&&&\\
&0.4&81.3&60.4&35.2&15.9&5.4&&&&&&99.8&95.7&72.9&29.2&4.6&&&&&\\\cdashline{2-22}[.4pt/1pt]
&0.6&99.2&96.3&89.3&75.9&55.4&5.5&&&&&100&100&100&99.6&95.2&5.2&&&&\\
&0.8&99.6&98.8&95.1&86.3&70.2&13.0&5.3&&&&100&100&100&99.9&98.4&15.3&4.8&&&\\
&1&99.9&99.7&99.2&96.7&89.7&32.9&19.6&5.1&&&100&100&100&100&99.9&61.8&40.0&5.3&&\\
&1.2&100&100&99.9&99.4&98.4&66.8&49.9&22.5&4.6&&100&100&100&100&100&96.1&89.7&44.4&5.1&\\
&1.4&100&100&100&99.9&99.9&92.3&83.5&56.4&34.1&5.1&100&100&100&100&100&99.9&99.7&90.9&51.4&4.7\\
\hline\hline
\end{tabular}
\caption{Estimated level ($d_1=d_2$) and power ($d_1\ne d_2$) (in \%) of the de-trended test based on
$\mathcal{I}_n$ with nominal level $\alpha = 5\%$. The samples are simulated from the model (\ref{changed})
with $\theta^* =1/2$. The sample sizes are $n=500$ and $n=5000$, the number of independent replications
is $10^4$.} \label{trend}
\end{sidewaystable}

\begin{sidewaystable}
      \centering \footnotesize
\begin{tabular}{l|l|lllll;{.4pt/1pt}lllll|lllll;{.4pt/1pt}lllll}
\hline\hline
&\multicolumn{1}{c}{}& \multicolumn{10}{c}{$n=500$} & \multicolumn{10}{c}{$n=5000$}\\
\cline{2-22}
$\underline\tau$ &\backslashbox{$d_2$}{$d_1$}&0&0.1&0.2&0.3&0.4&0.6&0.8&1&1.2&1.4&0&0.1&0.2&0.3&0.4&0.6&0.8&1&1.2&1.4\\
\hline\hline
0.05&0&2.3&&&&&&&&&&3.1&&&&&&&&&\\
&0.1&11.4&3.0&&&&&&&&&33.0&4.0&&&&&&&&\\
&0.2&30.1&12.2&3.7&&&&&&&&79.1&30.7&4.5&&&&&&&\\
&0.3&53.6&30.1&13.3&4.2&&&&&&&96.9&73.7&28.6&4.6&&&&&&\\
&0.4&71.9&49.6&28.3&13.0&4.8&&&&&&99.7&94.1&67.4&25.4&4.6&&&&&\\\cdashline{2-22}[.4pt/1pt]
&0.6&89.1&80.0&68.4&53.2&35.9&6.1&&&&&100&99.9&99.7&97.5&86.8&5.5&&&&\\
&0.8&91.1&83.9&73.6&59.8&44.6&7.5&5.8&&&&100&100&99.9&99.1&93.1&11.0&4.5&&&\\
&1&93.8&89.4&82.0&70.5&57.5&16.3&10.5&4.0&&&100&100&100&99.9&99.2&46.4&26.0&4.9&&\\
&1.2&97.3&94.6&90.6&83.0&74.5&34.1&24.3&14.5&4.0&&100&100&100&100&99.9&87.0&73.7&32.6&4.4&\\
&1.4&99.6&99.0&97.8&94.7&90.6&62.6&50.9&34.5&32.6&4.9&100&100&100&100&100&98.3&95.7&76.0&47.1&4.2\\\hline
0.1&0&2.7&&&&&&&&&&3.4&&&&&&&&&\\
&0.1&12.8&3.3&&&&&&&&&34.0&4.3&&&&&&&&\\
&0.2&32.7&13.5&4.0&&&&&&&&79.7&32.1&4.7&&&&&&&\\
&0.3&58.2&33.3&14.8&4.6&&&&&&&97.4&75.5&29.7&4.6&&&&&&\\
&0.4&78.0&55.8&32.4&14.5&4.9&&&&&&99.9&95.4&70.9&27.5&4.7&&&&&\\\cdashline{2-22}[.4pt/1pt]
&0.6&96.7&91.1&81.4&65.7&45.5&5.5&&&&&100&100&99.9&99.1&91.7&5.3&&&&\\
&0.8&97.8&94.7&88.1&75.6&58.6&8.8&5.3&&&&100&100&100&99.7&96.5&12.5&4.6&&&\\
&1&99.2&98.0&95.3&89.0&78.1&23.5&13.6&4.2&&&100&100&100&100&99.8&54.8&33.2&4.7&&\\
&1.2&99.8&99.5&98.7&96.9&92.2&51.6&36.5&18.5&4.1&&100&100&100&100&100&93.2&83.9&38.7&4.4&\\
&1.4&100&99.9&99.9&99.5&98.8&82.2&71.3&46.7&33.6&4.8&100&100&100&100&100&99.7&98.6&85.2&49.3&4.6\\\hline
0.2&0&3.2&&&&&&&&&&3.9&&&&&&&&&\\
&0.1&14.1&4.1&&&&&&&&&34.4&4.8&&&&&&&&\\
&0.2&35.5&15.0&4.8&&&&&&&&78.9&32.7&4.9&&&&&&&\\
&0.3&61.3&36.0&16.2&5.2&&&&&&&97.3&75.3&30.9&4.7&&&&&&\\
&0.4&81.9&59.6&35.6&15.7&5.3&&&&&&99.8&95.9&72.3&28.8&4.6&&&&&\\\cdashline{2-22}[.4pt/1pt]
&0.6&99.1&96.3&89.6&76.1&55.4&5.3&&&&&100&100&99.9&99.4&94.6&5.4&&&&\\
&0.8&99.6&98.6&95.1&86.2&70.8&12.5&5.2&&&&100&100&99.9&99.9&98.0&15.6&4.5&&&\\
&1&99.9&99.8&99.0&96.6&89.5&33.0&19.6&5.1&&&100&100&100&100&99.9&60.7&41.1&4.8&&\\
&1.2&100&99.9&99.9&99.7&98.4&67.1&49.2&22.8&4.6&&100&100&100&100&100&96.4&89.6&44.7&4.6&\\
&1.4&100&100&100&99.9&99.9&92.0&83.7&56.4&34.5&5.6&100&100&100&100&100&99.9&99.6&90.9&51.7&4.8\\
\hline\hline
\end{tabular}
\caption{Estimated level ($d_1=d_2$) and power ($d_1\ne d_2$) (in \%) of the de-trended
test based on $\mathcal{I}_n$ with nominal level $\alpha = 5\%$. The samples are simulated from
the model \eqref{eq:Xtrend} with $a=1$, $b=.01$ (for sample size $n=500$) and $a=1$, $b=.001$
(for sample size $n=5000$). The number of independent replications is $10^4$.} \label{trend2}
\end{sidewaystable}

\clearpage

\section{Appendix: proofs}\label{sec:app}

\noi {\it Proof of Proposition~\ref{general}.} (i) Without loss of generality, we will assume that
$A_{n1} = A_{n2} = 0$ in what follows. 
 Write
$z_{n1}(\tau) := \gamma^{-1}_{n1} S_{\lfloor n\tau\rfloor}$, $(n/\gamma_{n1}^2)
V_{\lfloor n\tau\rfloor}(X)  = \sum_{i=1}^6 U_{ni}(\tau)$, where the terms
\begin{eqnarray*}
U_{n1}(\tau)&:=&\big(n^2/\lfloor n\tau\rfloor^2\big) \int_0^{\lfloor n\tau\rfloor/n} z^2_{n1}(u) \d u, \\
U_{n2}(\tau)&:=&-2\big(n^2/\lfloor n\tau\rfloor^2\big) z_{n1}(\tau)
\int_0^{\lfloor n\tau\rfloor/n}  \big(\lfloor nu\rfloor/\lfloor n\tau\rfloor\big)
z_{n1}(u) \d u, \\
U_{n3}(\tau)&:=&\big(n^2/\lfloor n\tau\rfloor^2\big)  z^2_{n1}(\tau) \int_0^{\lfloor n\tau\rfloor/n}  \big(\lfloor nu\rfloor/\lfloor n\tau\rfloor\big)^2 \d u, \\
U_{n4}(\tau)&:=&
-\big(n^{3}/\lfloor n\tau\rfloor^{3}\big) \Big( \int_{0}^{\lfloor n\tau\rfloor/n} z_{n1}(u) \d u \Big)^2, \\
U_{n5}(\tau)&:=& 2\big(n^{3}/\lfloor n\tau\rfloor^{3}\big) z_{n1}(\tau)  \Big(
\int_{0}^{\lfloor n\tau\rfloor/n}  z_{n1}(u)  \d u \Big)
\int_0^{\lfloor n\tau\rfloor/n} \big(\lfloor nu\rfloor/\lfloor n\tau\rfloor\big) \d u, \\
U_{n6}(\tau)&:=& -\big(n^{3}/\lfloor n\tau\rfloor^{3}\big) z^2_{n1}(\tau)
\Big(\int_0^{\lfloor n\tau\rfloor/n} \big(\lfloor nu\rfloor/\lfloor n\tau\rfloor\big) \d u\Big)^2
\end{eqnarray*}
tend in distribution, as $n \to \infty$, to the corresponding limit quantities
\begin{eqnarray*}
U_{1}(\tau)&:=&\tau^{-2} \int_0^{\tau} Z^2_{1}(u) \d u, \\
U_{2}(\tau)&:=&-2\tau^{-2} Z_{1}(\tau) \int_0^{\tau}  (u/\tau)
Z_{1}(u) \d u, \\
U_{3}(\tau)&:=&\tau^{-2}  Z^2_{1}(\tau) \int_0^{\tau}  (u/\tau)^2 \d u, \\
U_{4}(\tau)&:=&
-\tau^{-3} \bigg(\int_{0}^{\tau} Z_{1}(u) \d u \bigg)^2, \\
U_{5}(\tau)&:=& 2\tau^{-3} Z_{1}(\tau)\bigg(\int_0^\tau Z_{1}(u)  \d u \bigg)
\int_0^{\tau} (u/\tau) \d u, \\
U_{6}(\tau)&:=& -\tau^{-3} Z^2_{1}(\tau) \bigg(\int_0^{\tau}
(u/\tau) \d u\bigg)^2.
\end{eqnarray*}
Note $Q_{\tau} (Z_1) = \sum_{i=1}^6 U_i(\tau)$ a.s. for each $\tau
\in (0, \upsilon_1]$. The joint convergence
\begin{equation} \label{Uconv}
(U_{n1}(\tau), \dots, U_{n6}(\tau)) \longrightarrow_{d}
(U_{1}(\tau), \dots, U_{6}(\tau))
\end{equation}
at each fixed point  $\tau \in (0, \upsilon_1]$ can be easily
derived from the (marginal) convergence $ \gamma^{-1}_{n1}
S_{\lfloor n\tau\rfloor}\, \longrightarrow_{D[0, \upsilon_1]}\, Z_1(\tau) $ in
(\ref{Dconvpower}). The convergence in (\ref{Uconv}) easily extends to
the joint convergence
at any finite number of points $0 < \tau_1 < \dots < \tau_m \le
\upsilon_1$. In other words,
\begin{equation} \label{V1}
(n/\gamma_{n1}^2) V_{\lfloor n\tau\rfloor}(X)
\longrightarrow_{{\rm fdd}(0,\upsilon_1]} Q_{\tau} (Z_1).
\end{equation}
In a similar way,
$$
\gamma^{-1}_{n2} S^*_{\lfloor n\tau\rfloor}\, \longrightarrow_{D[0,1-\upsilon_0]}\,
Z_2(\tau),
$$
implies
\begin{equation}\label{V2}
(n/\gamma_{n2}^2) V^*_{n-\lfloor n\tau\rfloor}(X)  \longrightarrow_{{\rm
fdd}[\upsilon_0,1)} Q_{1-\tau} (Z_2).
\end{equation}
It is clear from the joint convergence in (\ref{Dconvpower}) that
(\ref{V1}), (\ref{V2}) extend to the joint convergence of
finite-dimensional distributions, in other words, that
(\ref{VVlim}) holds with $\longrightarrow_{D(0, \upsilon_1] \times
D[\upsilon_0, 1)}$ replaced by $\longrightarrow_{{\rm fdd}(0,\upsilon_1] \times [\upsilon_0, 1)}$.

It remains to prove the tightness in $D(0, \upsilon_1] \times
D[\upsilon_0, 1)$. To this end, it suffices to check the
tightness of the marginal processes in (\ref{V1}) and (\ref{V2})
in the corresponding Skorokhod spaces $D(0, \upsilon_1]$ and
$D[\upsilon_0, 1)$. See, e.g., \cite{ferger}, \cite{whitt}.

Let us prove the tightness of the l.h.s. in (\ref{V1}) in $D(0,
\upsilon_1]$, or, equivalently, the tightness in $D[\upsilon,
\upsilon_1]$, for any $0< \upsilon < \upsilon_1$.
Let $\Upsilon_n(\tau) := (n/\gamma_{n1}^2) V_{\lfloor n\tau\rfloor}(X)$. Since $\{\Upsilon_n(\upsilon), n \ge 1\}$
is tight by (\ref{V1}), it suffices to show that for any $
\epsilon_1, \epsilon_2 >0$ there exist $\delta >0 $ and $n_0 \ge 1
$ such that
\begin{equation}\label{Qtight}
\P ( \omega_\delta (\Upsilon_n) \ge \epsilon_1 ) \le \epsilon_2, \qquad n
\ge n_0,
\end{equation}
where
$$
\omega_\delta(x) := \sup \big\{|x(a)-x(b)|: \,    \upsilon \le a <
b \le \upsilon_1, a-b < \delta \big\}
$$
is the continuity modulus of a function $x \in D[\upsilon,
\upsilon_1]$; see Billingsley (\citeyear{Billingsley:1968:CPM}, Theorem 8.2). Since
$\Upsilon_n(\tau) = \sum_{i=1}^6 U_{ni}(\tau)$, it suffices to show
(\ref{Qtight}) with $\Upsilon_n$ replaced by $U_{ni}, i=1, \dots, 6$,
in other words,
\begin{equation}\label{Utight}
\P ( \omega_\delta (U_{ni}) \ge \epsilon_1 ) \le \epsilon_2,
\qquad n \ge n_0,\quad i=1,\dots,6.
\end{equation}
We verify (\ref{Utight}) for $i=2$ only since the remaining cases
follow similarly. Write $U_{n2}(\tau) = \prod_{i=1}^3 H_{ni}(\tau)$, where
$H_{n1}(\tau) := -2\big(n^2/\lfloor n\tau\rfloor^2\big), \, H_{n2}(\tau) :=
z_{n1}(\tau), \, H_{n3}(\tau) := \int_0^{\lfloor n\tau\rfloor/n}
\big(\lfloor nu\rfloor/\lfloor n\tau\rfloor\big) z_{n1}(u) \d u$. Then
$\P (\omega_\delta (U_{n2}) \ge \epsilon_1 )
\le \sum_{i=1}^3 \big[\P ( \omega_\delta \big(H_{ni}) \ge
\epsilon_1/(3K)\big) + \P \big(\prod_{j\ne i} \|H_{nj}\| > K\big)
\big]$, where $\|x \| := \sup \{|x(a)|: \, \upsilon \le  a \le
\upsilon_1 \big\}$ is the sup-norm. Relation (\ref{Dconvpower})
implies that the probability $\P \big(\sum_{i=1}^3 \|H_{ni}\| > K\big)$
can be made arbitrary small for all $n > n_0(K)$ by a suitable
choice of $K$. By same relation (\ref{Dconvpower}) assumed under the
uniform topology,  for a given $\epsilon_1/K$,  we have that
$\lim_{\delta \to 0} \limsup_{n \to \infty} \P \big( \omega_\delta
(H_{ni}) \ge \epsilon_1/K\big) = 0$. This proves
(\ref{Utight}) and the functional convergence $(n/\gamma_{n1}^2)
V_{\lfloor n\tau\rfloor}(X) \longrightarrow_{D(0, \upsilon_1]} Q_{\tau} (Z_1)$.
The proof of $(n/\gamma_{n2}^2) V^*_{n-\lfloor n\tau\rfloor}(X)
\longrightarrow_{D[\upsilon_0, 1)} Q_{1-\tau} (Z_2)$ is
analogous. This concludes the proof of part (i), since the
continuity of the limit  process in (\ref{VVlim}) is immediate
from continuity of $(Z_1(\tau_1), Z_2(\tau_2)\big)$ and the
definition of $Q_\tau$ in (\ref{Q}).

\smallskip

\noindent (ii) Note that (\ref{bddbelow1}) and the a.s.\ continuity
of $\tau \mapsto Q_\tau (Z_1)$ guarantees that $\inf_{\tau \in
{\cal T}} Q_\tau (Z_1) >0$ a.s. Therefore relations (\ref{WZZconv})
follow from (\ref{VVlim}) and the continuous mapping theorem.
\smallskip

\noi (iii) Follows from (\ref{VVlim}) and the fact that $Z_1(\tau)
= 0, \, \tau \in {\cal T}$ implies $Q_\tau (Z_1) = 0, \, \tau \in
{\cal T}$. \hfill$\Box$

\bigskip

\noi {\it Proof of (\ref{approxF}).} Let first $-.5<d<.5$ and $b_i:=(a \star \pi(d))_i - \kappa \pi_i(d)$,
$i=0,1,\dots$. Consider the stationary process $\widetilde X_j := X_j -X^\dagger_j= \sum_{i=0}^\infty b_i\zeta_{j-i}$
with spectral density $\tilde f(x) = |\hat a(x)-\kappa|^2 g(x)$, where $\hat a(x)=\sum_{j=0}^\infty a_j \e^{-\i jx}$,
$\i := \sqrt{-1}$ and $g(x) :=(2\pi)^{-1} |1- \e^{-\i x}|^{-2d}$ is the spectral density of FARIMA$(0,d,0)$. We have
\begin{equation} \label{approxF1}
\E \bigg(\sum_{j=1}^n (X_j -X^\dagger_j) \bigg)^2 =  \int_{-\pi}^\pi \tilde f(x) D^2_n(x) \d x, \qquad D_n(x) := \frac{\sin (nx/2)}{\sin (x/2)}.
\end{equation}
Since $\hat a(x)$ is bounded and continuous on $[-\pi, \pi]$, $\hat a(0)=\kappa$, it follows that
$\tilde f(x) = o(|x|^{-2d})\, (x \to 0)$, which in turn implies (\ref{approxF}) for $-.5<d<.5 $;
see e.g.\ \citet[proof of Proposition~3.3.1]{giraitis3}.

Next, let $.5<d<1.5$. Then $X_j-X^\dagger_j=\sum_{k=1}^j \widetilde X_k$, where the stationary process
$\widetilde X_k:=\sum_{i=0}^\infty ((a \star \pi(d-1))_i-\kappa \pi_i(d-1)) \zeta_{k-i}$
satisfies $\E \big(\sum_{k=1}^j \widetilde X_k \big)^2 \le \epsilon(j) j^{2d-1}$, $\epsilon (j) \to 0$, $(j \to \infty)$,
see above. We have
\begin{eqnarray*}
 \E \bigg(\sum_{j=1}^n (X_j -X^\dagger_j)\bigg)^2 &=& \E \bigg(\sum_{j=1}^n \sum_{k=1}^j \widetilde X_k \bigg)^2
 \le \sum_{j_1,j_2=1}^n \bigg\{\E \bigg(\sum_{k=1}^{j_1} \widetilde X_k \bigg)^2  \E \bigg(\sum_{k=1}^{j_2} \widetilde X_k \bigg)^2 \bigg\}^{1/2}\\
 &\le&  \bigg\{\sum_{j=1}^n \sqrt{\epsilon(j) j^{2d-1}} \bigg\}^2 =  o(n^{2d+1}).
\end{eqnarray*}
This completes the proof of (\ref{approxF}).\hfill $\Box$
\bigskip

\noi {\it Proof of Proposition~\ref{linearII}.} Note first that the convergence
in (\ref{DconvII}) for type II integrated process of (\ref{IId}) can be easily established
following the proof of Proposition~\ref{linear}. Hence, it suffices to show that
\begin{equation}\label{initial00}
 \sup_{\tau \in [0,1]} \frac{1}{n^{d+.5}}\bigg|\sum_{t=1}^{\lfloor n\tau\rfloor} R^0_t\bigg|\ \rightarrow_p \  0,
\end{equation}
where
$R^0_t := (1-L)_+^{-d} (1-L)_-^d X^0_t = \sum_{i=0}^\infty X^0_{-i} \sum_{j=0}^{t-1} \pi_j(d) \pi_{t-j+i}(-d)$
is the contribution arising from initial values. When  $d>0.5$, using (\ref{initial}), the Cauchy-Schwarz
inequality and the fact that $|\pi_j(d)|\le Cj^{d-1}$ we obtain
\begin{eqnarray*}
\E \bigg(\frac{1}{n^{d + .5}}\sum_{t=1}^n |R^0_t|\bigg)^2
&\le&\frac{C}{n^{2d+1}} \bigg(\sum_{t=1}^n \sum_{i=0}^\infty \sum_{j=0}^{t-1} |\pi_j(d)| \, |\pi_{t-j+i}(-d)|\bigg)^2 \nonumber \\
&\le&\frac{C}{n^{2d+1}} \bigg(\sum_{1 \le j < t \le n} j^{d-1} (t-j)^{-d}\bigg)^2  \ \to \  0 \label{initial0}
\end{eqnarray*}
since
\begin{eqnarray*}
\sum_{1 \le j < t \le n} j^{d-1} (t-j)^{-d} \le \sum_{j=1}^n j^{d-1} \sum_{k=1}^n k^{-d} \  \le \
C\begin{cases}
n^d, &d> 1, \\
n \log n, &d = 1, \\
n, &\frac 12< d< 1.
\end{cases}
\end{eqnarray*}
This proves (\ref{initial00}). \hfill $\Box$

\bigskip

\noi {\it Proof of Proposition~\ref{power}.} We restrict the  proof to the case (i) and $i=2$,
or, equivalently, to the model (\ref{X2t}), since the remaining cases 
can be treated similarly. Similarly as in the proof of (\ref{VVlim}), it suffices to prove the joint convergence of finite-dimensional
distributions in (\ref{ZZ1}) and the functional convergence of marginal processes, viz.,
\begin{eqnarray}\label{ZZiii}
n^{-d_1-.5} S_{2,\lfloor n\tau\rfloor}&\longrightarrow_{D[0,\theta]}&Z_{2,1}(\tau),  \qquad
n^{-d_2-.5} S^*_{2,\lfloor n\tau\rfloor}\ \longrightarrow_{D[0, 1-\underline\tau]}\ Z_{2,2}(\tau).
\end{eqnarray}
Since $X_{2,t}  = \sum_{j=0}^t \pi_j(d_1) \zeta_{t-j}, \, 1\le t \le \lfloor n\underline\tau\rfloor$ has constant memory parameter $d_1$,
the proof of the first convergence in (\ref{ZZiii}) to $Z_{2,1}(\tau) = B^{\rm II}_{d_1+.5}(\tau) $ is standard,
 and we omit it. 

Consider the second convergence in (\ref{ZZiii}). It can be rewritten as
\begin{eqnarray}\label{ZZ4}
&n^{-d_2-.5} S_{2,\lfloor n\tau\rfloor} \,
\longrightarrow_{D[\underline\tau,1]}\, {\cal Z}_2(\tau),
\end{eqnarray}
where $S_{2,\lfloor n\tau\rfloor}=\sum_{t=1}^{\lfloor n\tau\rfloor} X_{2,t} = \sum_{t=1}^{\lfloor n\tau\rfloor} \sum_{s=1}^t b_{2,t-s}(t) \zeta_s.$
Let us first prove the one-dimensional convergence in (\ref{ZZ4}) at a fixed point $\tau \in [\underline\tau, 1]$.

We start with the case $\tau>\overline\theta$. Following the scheme of discrete stochastic
integrals in \cite{surg1}, rewrite the l.h.s.\ of (\ref{ZZ4}) as a discrete stochastic integral
$$
n^{-d_2-.5} S_{2,\lfloor n\tau\rfloor} = \int_0^\tau F_n(u) \d z_n(u) = \int_0^{\overline\theta} F_n(u) \d z_n(u)
+ \int_{\overline\theta}^\tau  F_n(u) \d z_n(u),
$$
where $z_n(u):= n^{-1/2} \sum_{i=1}^{\lfloor nu\rfloor}\zeta_i$ is the partial sum process of standardized i.i.d.\ r.v.s,
tending weakly to a Brownian motion $\{B(u), u\in [0,1]\}$. The integrand $F_n$ in the above integral is equal to
\begin{eqnarray*}
F_n(u)&:=&n^{-d_2}\sum_{t= \lfloor nu\rfloor}^{\lfloor n\tau\rfloor} b_{2,t-\lfloor nu\rfloor}(t) \\
&=&
\begin{cases}n^{-d_2}\sum_{t= \lfloor nu\rfloor}^{\lfloor n\tau\rfloor} b_{2,t-\lfloor nu\rfloor}(t), &0< u \le \overline\theta, \\
n^{-d_2}\sum_{t=\lfloor nu\rfloor}^{\lfloor n\tau\rfloor} \pi_{t-\lfloor nu\rfloor}(d_2), &\overline\theta < u \le \tau,
\end{cases}
\end{eqnarray*}
where we used the fact that $b_{2,t-\lfloor nu\rfloor}(t) = \pi_{t-\lfloor nu\rfloor}(d_2)$ for $t \ge \lfloor nu\rfloor \ge \lfloor
n\overline \theta\rfloor$, where $\pi_j(d)$ are the FARIMA coefficients in (\ref{farima}).
Similarly, the r.h.s.\ of (\ref{ZZ4}) can be written as the sum of two stochastic integrals:
$$
\int_0^\tau F(u) \d B(u) = \int_0^{\overline\theta} F(u) \d B(u) + \int_{\overline\theta}^\tau  F(u) \d B(u),
$$
where
\begin{eqnarray*}
F(u)&:=&
\begin{cases}\Gamma(d_2)^{-1} \int_{\overline\theta}^\tau (v-u)^{d_2-1} \e^{H(u,v)} \d v, &0< u \le \overline\theta, \\
\Gamma(d_2+1)^{-1}(\tau - u)^{d_2}, &\overline\theta < u \le \tau.
\end{cases}
\end{eqnarray*}

Accordingly, using the above mentioned criterion in Surgailis (\citeyear{surg1}, Proposition~3.2)
(see also Lemma~2.1 in \cite{bruzaite}), the one-dimensional convergence
in (\ref{ZZ4}) follows from the $L^2-$con\-ver\-gence of the integrands:
\begin{equation} \label{2conv}
\int_0^{\overline\theta} |F_n(u) - F(u)|^2 \d u \to 0, \qquad \int_{\overline\theta}^\tau |F_n(u) - F(u)|^2 \d u \to 0.
\end{equation}
The second relation in (\ref{2conv}) is easy using the properties of FARIMA filters. Denote $J_n$ the first integral
in  (\ref{2conv}). The integrand there can be rewritten as
\begin{equation}\label{decomposition}F_n (u) - F(u) = \int_u^{\overline\theta} n^{1-d_2} b_{2,\lfloor nv\rfloor-\lfloor nu\rfloor}(\lfloor nv\rfloor) \d v \, +\, \int_{\overline\theta}^\tau G_n(u,v) \d v   \, -\,  n^{-d_2}(n\tau - \lfloor n\tau\rfloor)   b_{2,\lfloor n\tau\rfloor-\lfloor nu\rfloor}(\lfloor n\tau\rfloor),\end{equation}
 where $G_n(u,v):=n^{1-d_2} b_{2,\lfloor nv\rfloor-\lfloor nu\rfloor}(\lfloor nv\rfloor) -\Gamma(d_2)^{-1} (v-u)^{d_2-1} \e^{H(u,v)}$ and $H(u,v)$ is defined at (\ref{Huv}).

Let us write $b_{2,\lfloor nv\rfloor-\lfloor nu\rfloor}(\lfloor nv\rfloor)= \pi_{\lfloor nv\rfloor-\lfloor nu\rfloor}(d_2) K_n(u,v)$ where
\begin{eqnarray} \label{Kuv}
K_n(u,v) := \frac{b_{2,\lfloor nv\rfloor-\lfloor nu\rfloor}(\lfloor nv\rfloor)}{\pi_{\lfloor nv\rfloor-\lfloor nu\rfloor}(d_2)}
\ = \
\prod_{i = 1 \vee (\lfloor nv\rfloor-\lfloor n \overline\theta\rfloor) }^{\lfloor nv\rfloor - \lfloor nu\rfloor} \frac{d\big(\frac{\lfloor nv\rfloor-i+1}{n}\big) - 1 + i}{d_2 - 1 + i}.
\end{eqnarray}
We claim that
\begin{equation}\label{Klim}
\lim_{n \to \infty} K_n(u,v)= \e^{ H(u,v)}, \qquad 0\le u \le \overline\theta <  v \le 1,
\end{equation}
Indeed,
\begin{eqnarray*}
K_n(u,v)
&=&\exp \Bigg\{ \sum_{i = 1 \vee (\lfloor nv\rfloor-\lfloor n\overline\theta\rfloor)}^{\lfloor nv\rfloor - \lfloor nu\rfloor}
\log \bigg(1 -  \frac{d_2 - d\big(\frac{\lfloor nv\rfloor-i}{n}\big) }{d_2 - 1 + i}\bigg) \Bigg\} \
=\  \e^{H_n(u,v) + R_n(u,v)},
\end{eqnarray*}
where
\begin{eqnarray*}
H_n(u,v)&:=&n^{-1}\sum_{i = 1\vee (\lfloor nv\rfloor-\lfloor n\overline\theta\rfloor)}^{\lfloor nv\rfloor - \lfloor nu\rfloor} \frac{d\big(\frac{\lfloor nv\rfloor-i}{n}\big)-d_2}
{\frac{d_2 - 1 + i}{n}} \ \to \ H(u,v), \\
R_n(u,v)&=& O\Bigg(\sum_{i = 1\vee (\lfloor nv\rfloor-\lfloor n\overline\theta\rfloor)}^{\lfloor nv\rfloor - \lfloor nu\rfloor} \frac{1}{i^2}\Bigg) \ = \  O\bigg(\frac{1}{1 \vee (\lfloor nv\rfloor-\lfloor n\overline\theta\rfloor)}\bigg),
\end{eqnarray*}
hence $R_n(u,v) \to 0$ for any $v > \overline\theta$.

 The proof of $J_n \to 0$ in (\ref{2conv}) then follows from the following arguments. Using on one hand the fact that the ratio $K_n(u,v)$  tends to 0 for $0<u< v \le \overline\theta$, on the other hand (\ref{Klim}), and from the well-known asymptotics  $\pi_j(d) \sim \Gamma(d)^{-1} j^{d-1}, \, j\to \infty$ of FARIMA coefficients, it easily follows that $n^{1-d_2} b_{2,\lfloor nv\rfloor-\lfloor nu\rfloor}(\lfloor nv\rfloor)\to 0$ for any  $0<u< v \le \overline\theta$, and  $G_n(u,v) \to 0 $ for any $0< u  < v \le 1 $ fixed.  Moreover, the last term in \eqref{decomposition} obviously tends to 0 because $d_2>0$. Since
both sides of (\ref{Klim}) are nonnegative and bounded by 1, the above convergences extend to
the proof of $J_n \to 0$ by the dominated convergence theorem. This proves the convergence
of one-dimensional distributions in (\ref{ZZ4}) for $\tau > \overline\theta$.

For $\underline\tau \le \tau \le \overline\theta$, the above convergence
follows similarly by using the fact that  $K_n(u,v) $  tends to 0 for $0<u< v \le \overline\theta$.

The proof of the convergence of general finite-dimensional distributions in (\ref{ZZ4}), as well as
 the joint convergence of finite-dimensional distributions in (\ref{ZZ1}),
can be achieved analogously, by using the Cram\`er-Wold device. Finally, the tightness
in (\ref{ZZ4}) follows by the Kolmorogov criterion (see, e.g. \cite{bruzaite}, proof of Theorem 1.2 for details).
Proposition~\ref{power} is proved.\hfill $\Box$

\bigskip

\noi {\it Proof of Proposition~\ref{VVtrend}.}
Consider de-trended observations and their partial sums processes as defined by
\begin{eqnarray*}
\veps_j&:=&X_j - (A_{n}- B_{n})  - 2 B_{n} j, \qquad
S_k(\veps) \ :=\  \sum_{j=1}^k \veps_j, \qquad S^*_{n-k}(\veps) \ :=\  \sum_{j=k+1}^n \veps_j.
\end{eqnarray*}
Note that $S_k(\veps) = S_k(X) - k (A_{n} -B_{n})   - k (k+1) B_{n}
= S_k(X) - k A_{n} - k^2 B_{n}$ and the null hypothesis  ${\bf H}^{\rm trend}_0$ can be rewritten as
\begin{eqnarray}\label{Dvepsconv1}
 \gamma^{-1}_{n} S_{\lfloor n\tau\rfloor}(\veps)
 &\longrightarrow_{D[0,1]}& Z(\tau).
\end{eqnarray}
For a fixed $1\le k < n$, let  $(\hat a_k(\veps), \hat b_k(\veps)) := {\rm argmin} \big(\sum_{j=1}^k (\veps_j - a - bj)^2 \big), \
(\hat a^*_{n-k}(\veps), \hat b^*_{n-k}(\veps)) \ := \ {\rm argmin} \big(\sum_{j=k+1}^n (\veps_j - a - bj)^2\big)$
be the corresponding linear regression coefficients. More explicitly,
\begin{eqnarray}\label{hatab}
\hat a_k(\veps)&=&\frac{\Big(\sum\limits_{j=1}^k j\Big)\Big(\sum\limits_{j=1}^k j \veps_j\Big) -
\Big(\sum\limits_{j=1}^k j^2\Big) \Big(\sum\limits_{j=1}^k  \veps_j\Big)}
{\Big(\sum\limits_{j=1}^k j\Big)^2 - k \sum\limits_{j=1}^k j^2 }, \qquad
\hat b_k(\veps)
\ :=\  \frac{\Big(\sum\limits_{j=1}^k j\Big)\Big(\sum\limits_{j=1}^k \veps_j\Big) -
k \sum\limits_{j=1}^k j \veps_j}
{\Big(\sum\limits_{j=1}^k j\Big)^2 - k \sum\limits_{j=1}^k j^2 }.
\end{eqnarray}
It is easy to verify that $\hat a_k(\veps)
= \hat a_k(X) + (B_{n} - A_{n}), \
\hat b_{k}(\veps)\ =\  \hat b_k(X) - 2B_{n}$. Hence
we obtain the following expression of residual partial sums $\widehat S_j = S_j(\widehat X)$ via de-trended partial sums $S_j(\veps)$
and the above regression coefficients:
\begin{eqnarray}\label{hatS}
\widehat S_j
&=&S_j(\veps) - j \Big(\hat a_k(\veps) + \frac{\hat b_k(\veps)}{2}\Big)  - j^2 \Big(\frac{\hat b_k(\veps)}{2}\Big)
\end{eqnarray}
The limit behavior of  ${\mathcal V}_{\lfloor n\tau\rfloor}(X)$ follows from the limit behavior
of $ \gamma_n^{-1}\widehat S_{\lfloor nu\rfloor}, u \in [0,\tau]$ similarly as in the proof of
Proposition~\ref{general}. The behavior of the first term  $\gamma_n^{-1}S_{\lfloor nu\rfloor}(\veps)$
in (\ref{hatS}) is given in (\ref{Dvepsconv1}). It remains to identify the limit regression
coefficients $\hat a_{\lfloor n\tau\rfloor}(\veps)$,  $\hat b_{\lfloor n\tau\rfloor}(\veps)$ in (\ref{hatS}).
Clearly the denominator $\big(\sum_{j=1}^k j\big)^2 - k \sum_{j=1}^k j^2 \sim -\frac{k^4}{12}$.
The numerators in (\ref{hatab}) are written in terms of $S_k(\veps)$ and $\sum_{j=1}^k j \veps_j$.
From summation by parts and (\ref{Dvepsconv1}) we obtain
\begin{eqnarray}\label{Dvepsconv2}
n^{-1} \gamma_n^{-1} \sum_{j=1}^{\lfloor n\tau\rfloor} j \veps_j
&=&n^{-1} \gamma_n^{-1}\bigg( \lfloor n\tau\rfloor S_{\lfloor n\tau\rfloor}(\veps) - \sum_{j=1}^{\lfloor n\tau\rfloor-1} S_j (\veps)\bigg) \
\longrightarrow_{D[0,1]} \  \tau Z(\tau) - \int_0^\tau Z(v) \d v.
\end{eqnarray}
Relations (\ref{Dvepsconv1}), (\ref{hatab})  and (\ref{Dvepsconv2}) entail $(\lfloor nu\rfloor/\gamma_n) \hat a_{\lfloor n\tau\rfloor}(\veps)
\rightarrow_{D[0,\tau]}u \hat a_\tau$, $(\lfloor nu\rfloor^2/(2\gamma_n))\hat b_{\lfloor n\tau\rfloor}(\veps)
\rightarrow_{D[0,\tau]} u^2 \frac{\hat b_\tau}{2}$, where
\begin{eqnarray}\label{ablim}
\hat a_\tau \ := \  - 2\Big(\frac{1}{\tau}\Big) Z(\tau) + 6 \Big(\frac{1}{\tau}\Big)^2  \int_0^\tau Z(v) \d v, \qquad
\hat b_\tau \ := \
6\Big(\frac{1}{\tau}\Big)^2 Z(\tau) - 12 \Big(\frac{1}{\tau}\Big)^3  \int_0^\tau Z(v) \d v,
\end{eqnarray}
leading to the convergence $\gamma_n^{-1}\widehat S_{\lfloor nu\rfloor}\longrightarrow_{D[0,\tau]}
{\mathcal Z}(u,\tau) $, where the limit process is given in (\ref{Qtrend}). The remaining details of the proof
are similar as in Proposition~\ref{general}.
\hfill$\Box$

\newpage

\bibliographystyle{apalike} 



\end{document}